\documentclass[a4paper,12pt]{article}

\usepackage[hyperref,english]{package_RFb}
\usepackage{numbered_constants}
\usepackage[page,title]{appendix}

\usepackage{geometry}
\usepackage[numbers]{natbib}
\title{Stochastic partial differential equations describing neutral genetic diversity under short range and long range dispersal}
\author{Rapha\"el Forien\footnote{INRAE, BioSP, 84914, Avignon, France, e-mail: raphael.forien@inrae.fr}}

\RenewDocumentCommand{\1}{}{\mathds{1}}

\DeclareConstant{Cst}{C}

\begin{document}
	\maketitle
	
	\begin{abstract}
		In this paper, we consider a mathematical model for the evolution of neutral genetic diversity in a spatial continuum including mutations, genetic drift and either short range or long range dispersal.
		The model we consider is the spatial $ \Lambda $-Fleming-Viot process introduced by Barton, Etheridge and Véber, which describes the state of the population at any time by a measure on $ \R^d \times [0,1] $, where $ \R^d $ is the geographical space and $ [0,1] $ is the space of genetic types.
		In both cases (short range and long range dispersal), we prove a functional central limit theorem for the process as the population density becomes large and under some space-time rescaling.
		We then deduce from these two central limit theorems a formula for the asymptotic probability of identity of two individuals picked at random from two given spatial locations.
		In the case of short range dispersal, we recover the classical Wright-Malécot formula, which is widely used in demographic inference for spatially structured populations.
		In the case of long range dispersal we obtain a new formula which could open the way for a better appraisal of long range dispersal in inference methods.
	\end{abstract}
	
	 \section*{Introduction}
	 
	 \paragraph*{Isolation by distance patterns}
	 
	 Many populations occupy a geographical area whose spatial extent is much larger than the typical distance travelled by individuals during their lifetime.
	 As a result, individuals living close to each other are on average more related than those living far apart.
	 Genetic similarity between individuals thus decreases as a function of the geographical distance between them.
	 This is known as \textit{isolation by distance}, and the exact speed and shape of this decrease of genetic similarity has been the subject of many theoretical studies (\cite{wright_isolation_1943,wright_breeding_1940,malecot_heterozygosity_1975}, followed by \cite{kimura_stepping-stone_1953,kimura_stepping_1964,sawyer_continuous_1981} and many others).
	 
	 The quantity that has been most often used to describe this phenomenon is the \textit{probability of identity by descent}, which is defined as follows.
	 Given a model describing the evolution of the genetic composition of a spatially structured population, the probability of identity by descent is the probability that two individuals sampled from two distinct locations carry the same genetic material (or allele) at a given locus and that this allele was inherited from a common ancestor without any mutation between them.
	 If the dispersion of individuals in the population is very local, this probability is approximately given by the so-called Wright-Malécot formula, which takes the following form \cite{barton_neutral_2002}.
	 If $ x > 0 $ is the distance between the two sampled individuals, $ \mu $ is the mutation rate of individuals and $ \sigma^2 $ is the average square distance between an individual and its parent, then, in a two dimensional space,
	 \begin{equation} \label{original_wright_malecot}
	 	\frac{1}{\mathcal{N} + \log(\sigma/(\kappa\sqrt{2\mu}))} K_0\left( \sqrt{2\mu} \frac{x}{\sigma} \right), \quad  x > \kappa,
	 \end{equation}
	 where $ \kappa $ is a parameter called the local scale and $ K_0(\cdot) $ is the modified Bessel function of the second kind of degree 0.
	 On the other hand, in a one dimensional space,
	 \begin{align*}
	 \frac{1}{1+\frac{\mathcal{N} \sqrt{2\mu}}{\pi \sigma}} \exp \left( - \sqrt{2\mu} \frac{x}{\sigma} \right).
	 \end{align*}
	 The constant $ \mathcal{N} $ is called Wright's neighbourhood size and is defined by
	 \begin{align*}
	 	\mathcal{N} = 2\pi \sigma^2 N,
	 \end{align*}
	 where $ N $ is proportional to the density of breeding individuals in the population.
	 The parameter $ \kappa $ depends on the details of the specific model under consideration, and is chosen so that the value of \eqref{original_wright_malecot} at $ x = \kappa $ matches the probability of identity by descent of a pair of individuals sampled close to each other (see \cite{barton_neutral_2002} for more details).
	 A formula similar to \eqref{original_wright_malecot} appears in the correlations between allele frequencies at two spatial locations \cite{kimura_stepping_1964}, as expected from the duality relation between genealogical and forwards-in-time models \cite{barton_modeling_2013}.
	 
	 The Wright-Malécot formula allows one - if a large enough number of individuals are sequenced at a sufficient number of loci - to estimate the demographic parameters $ \mathcal{N} $ and $ \sigma $ of real populations by fitting \eqref{original_wright_malecot} to pairwise identity matrices \cite{rousset_genetic_1997,barton_inference_2013}.
	 More sophisticated inference methods have also been recently developed, using long continuous tracts of shared genetic material (called blocks of identity by descent, or IBD blocks), and these also rely on the Wright-Malécot formula \cite{ringbauer_inferring_2017}.
	 
	 In some species, however, individuals can disperse their offspring arbitrarily far away from their own location, and the average square distance between an offspring and its parent might be very large, violating one assumption of the Wright-Malécot formula.
	 Evidence of long range dispersal has been found for example in plant species \cite{cain_long-distance_2000} and fungi \cite{buschbom_migration_2007} (see also \cite{nathan_methods_2003}), but, up to now, no analogue formula was available for such populations.
	 One of the aims of this paper is to fill this gap.
	 
	 \paragraph*{Modelling evolution in a spatial continuum}
	 
	 In this paper, we use a model called the spatial $ \Lambda $-Fleming-Viot process (SLFV in short), introduced by N. Barton, A. Etheridge and A. Véber in \cite{etheridge_drift_2008,barton_new_2010-1} (see also \cite{barton_modeling_2013} for a review).
	 This model describes the genetic composition of a spatially structured population by a measure on $ \R^d \times [0,1] $, where $ \R^d $ is the geographical space and $ [0,1] $ is the space of genetic types.
	 The density of the population is tightly regulated so that, at any time, the spatial marginal of the measure is always the Lebesgue measure.
	 
	 The population then evolves through a sequence of reproduction events, where each of these events affects a ball drawn according to a Poisson point process.
	 During these events, a fixed proportion of the individuals in the ball die and are replaced by the offspring of an individual chosen at random within the same ball.
	 We also include mutations by assuming that individuals change their type to a new one, chosen uniformly in $ [0,1] $, at some fixed rate $ \mu > 0 $.
	 As a result of this construction, if we start the process from the Lebesgue measure on $ \R^d \times [0,1] $, two individuals share the same genetic type if and only if they share a common ancestor which is more recent than the last mutation to occur in their genealogy.
	 
	 This way of encoding the genetic composition of a spatially structured population as a measure on $ \R^d \times [0,1] $ was already used in several settings, \textit{e.g.} in \cite{handa_measure-valued_1990,evans_coalescing_1997,liang_two_2009}.
	 This particular way of including mutations in the SLFV was also introduced in \cite{veber_spatial_2015}.
	 In particular, the SLFV records the size and geographical extent of all the families in the population, see Remark~\ref{remark:identity_in_state} below.
	 
	 We consider two separate settings: one in which the radius of reproduction events is fixed, which corresponds to local dispersal, and one in which this radius is drawn at random according to some heavy-tailed distribution, corresponding to long range dispersal.
	 
	 \paragraph*{Large population - rare mutations limit}
	 
	 We then consider the limit of this process as both the mutation rate and the fraction of individuals replaced at each reproduction event converge to zero.
	 We show that, if we rescale time and space properly, the SLFV converges to the Lebesgue measure on $ \R^d \times [0,1] $.
	 In other words, the probability that any two given individuals are related vanishes.
	 To recover isolation by distance patterns, we need to look at the fluctuations of the process around its deterministic limit.
	 
	 We do this by proving a central limit theorem for the SLFV.
	 More precisely, we show that the rescaled difference between the SLFV and its limit converges to a distribution-valued Gaussian process, given as the solution to a linear stochastic partial differential equation on $ \R^d \times [0,1] $.
	 The coefficients of this SPDE are functions of the parameters of the SLFV, and more importantly, the driving noise depends on the dispersal assumption. In the case of a fixed radius of reproduction events (\textit{i.e.} short range dispersal, Theorem~\ref{thm:clt_fixed_radius} below), the driving noise is white in space and time, and has a Fleming-Viot component at each spatial location.
	 In the case of heavy-tailed radius of reproduction events (long range dispersal, Theorem~\ref{thm:clt_stable}), the driving noise becomes correlated in space and the strength of these correlations depends on the decay of the radius distribution.
	 
	 These results extend previous results obtained in \cite{forien_central_2017}, where rescaling limits as well as central limit theorems were obtained for the two-types SLFV with natural selection.
	 The main difficulty in the present paper compared to \cite{forien_central_2017} is the fact that, while the SLFV with selection took values in a space  of measurable maps from $ \R^d $ to $ [0,1] $, the SLFV with mutations takes values in a space of measurable maps from $ \R^d $ to $ \mathcal{M}_1([0,1]) $, the space of probability measures on $ [0,1] $ (see Section~\ref{subsec:state_space} below).
	 
	 \paragraph*{The Wright-Malécot formula}
	 
	 We are then able to use our results on the asymptotic fluctuations of the SLFV with mutations around the Lebesgue measure on $ \R^d \times [0,1] $ to deduce the asymptotic behaviour of the probability of identity of two individuals sampled from two different regions.
	 The formula for this probability is obtained by computing the singular part of the Gaussian process given by the central limit theorem for the SLFV with mutations.
	 In the case of short range dispersal, we recover the classical Wright-Malécot formula \eqref{original_wright_malecot}, albeit without the term involving $ \kappa $.
	 This is because we are considering a limit where the population density tends to infinity, hence this term becomes negligible compared to $ \mathcal{N} $ in our result.
	 In the case of long range dispersal, we obtain with the same method a new formula for the probability of identity of two individuals sampled from the population, as a function of the decay of the radius distribution.
	 
	 It is worth noting that, in this paper, we never use the genealogical dual of the SLFV, but we are still able to obtain estimates of the probability of identity of two individuals in the population.
	 This is important because, in some situations, a proper genealogical dual may be hard to find, or even non-existent (see for example \cite{biswas_spatial_2021}).
	 Our techniques might then allow one to recover genealogical information about certain populations even in the absence of a dual.
	 
	 \paragraph*{Demographic inference}
	 
	 We hope that this result will permit new developments in demographic inference methods in order to better take into account long range dispersal in natural populations.
	 Current parametric estimation methods could be adapted to include this effect and estimate the strength of long range dispersal as an additional parameter (through the exponent of the fractional Laplacian appearing in the limiting equations below).
	 
	 ~~
	 
	 The paper is laid out as follows.
	 In Section~\ref{sec:model}, we define the SLFV with mutations.
	 In Section~\ref{sec:results}, we state our main results, namely two central limit theorems (one for short range dispersal and one for long range dispersal) and we give the two corresponding formulas for the probability of identity.
	 We also state a central limit theorem for the SLFV with mutations in a slightly more general setting which is of independent interest.
	 The two central limit theorems are proved in Section~\ref{sec:proof_clt}.
	 We introduce general notations which allow us to prove the two results at the same time, following a general strategy which is outlined in Section~\ref{subsec:strategy}.
	 The two formulas for the probability of identity are then proved in Section~\ref{sec:WM_proof}.
	 Finally, in the Appendix, we recall a few useful results and we show how to adapt a particular result of \cite{walsh_introduction_1986} to our setting.
	 
	 \section{Definition of the model} \label{sec:model}
	 
	 Consider a population occupying a continuous geographical space (here $ \R^d $) and where each individual carries a genetic type belonging to $ [0,1] $.
	 The state of the population at time $ t \geq 0 $ can then be represented by a (random) map
	 \begin{align*}
	 	\rho_t : \R^d \to \mathcal{M}_1([0,1]),
	 \end{align*}
	 where $ \mathcal{M}_1([0,1]) $ denotes the space of probability measures on $ [0,1] $.
	 In this way, $ \rho_t(x, dk) $ is the probability measure corresponding to the distribution of genetic types at $ x \in \R^d $ at time $ t \geq 0 $.
	 More precisely, for $ A \subset [0,1] $, 
	 \begin{align*}
	 	\int_A \rho_t(x, dk)
	 \end{align*}
	 is the probability that an individual sampled uniformly from those present at $ x \in \R^d $ at time $ t $ carries a type belonging to $ A $.
	 The evolution of $ (\rho_t, t \geq 0) $ is then governed by births and deaths in the population, along with migration, mutations, and possibly natural selection.
	 Barton, Etheridge and Véber introduced the \slfv (SLFV in short) in \cite{barton_new_2010-1} as a framework to describe the evolution of $ (\rho_t, t \geq 0) $ as individuals die and reproduce in the population.
	 
	 \subsection{The state space of the process} \label{subsec:state_space}
	 
	 Before giving the formal definition of the SLFV, let us precise its state space and its topology.
	 In particular, we identify two measurable maps from $ \R^d $ to $ \mathcal{M}_1([0,1]) $ if they coincide up to a set of Lebesgue measure zero.
	 Let $ \Xi $ denote the corresponding quotient space.
	 Then $ \Xi $ is in one-to-one correspondence with the set of non-negative Radon measures on $ \R^d \times [0,1] $ whose spatial marginal is the Lebesgue measure on $ \R^d $ \cite[Section~2.2]{veber_spatial_2015} via the relation
	 \begin{align} \label{xi_measure_oto_correspondence}
	 	m(dx\, dk) = \rho(x, dk) dx.
	 \end{align}
	 For $ \rho \in \Xi $ and $ \phi : \R^d \times [0,1] \to \R $, let $ \langle \rho, \phi \rangle $ denote the integral of $ \phi $ against the measure $ \rho $.
	 
	 We equip the space $ \Xi $ with the topology of vague convergence and the associated Borel $ \sigma $-field (\textit{i.e.} $ \rho_{n} \to \rho \in \Xi $ if and only if $ \langle \rho_n, \phi \rangle \to \langle \rho, \phi \rangle $ for any $ \phi : \R^d \times [0,1] \to \R $ that is continuous and compactly supported).
	 Endowed with this topology, the space $ \Xi $ is both compact and metrisable \citep[Lemmas~2.1 and 2.2]{veber_spatial_2015}.
	 More precisely, if $ (\phi_n, n \geq 1) $ is a sequence of uniformly bounded continuous and compactly supported functions $ \phi_n : \R^d \times [0,1] \to \R $ which separates points in $ \Xi $ (in the sense that $ \langle \rho, \phi_n \rangle = \langle \rho', \phi_n \rangle $ for all $ n \geq 1 $ if and only if $ \rho = \rho' $), then
	 \begin{align} \label{def:d}
	 \quad d(\rho, \rho') = \sum_{n=1}^{\infty} \frac{1}{2^n} \abs{\langle \rho, \phi_n \rangle - \langle \rho', \phi_n \rangle}, \quad \forall \rho, \rho' \in \Xi,
	 \end{align}
	 defines a metric on $ \Xi $ which induces the vague topology.
	 
	 For $ \phi : \R^d \times [0,1] \to \R $ and $ q \geq 1 $, define
	 \begin{align} \label{def_norm_q}
	 \| \phi \|_q = \left( \int_{\R^d} \sup_{k \in [0,1]} \abs{\phi(x,k)}^q dx \right)^{1/q}
	 \end{align}
	 whenever the right hand side is finite, and let $ E_q $ be the space of all measurable real-valued functions $ \phi $ such that $ \| \phi \|_q < \infty $ and $ \sup_{(x, k) \in \R^d \times [0,1]} \abs{\phi(x,k)} < \infty $.
	 Also for a multi-index $ \beta \in \N^d $, let $ \partial_\beta \phi $ denote the partial derivative of $ \phi : \R^d \times [0,1] \to \R $ with respect to the space variable, and let $ \abs{\beta} = \sum_{i=1}^{d} \beta_i $.
	 We can (and do in the rest of the paper) assume that the $ \phi_n $ are all smooth and that there exists a constant $ \newCst{phi_n} > 0 $ such that, for all $ n \geq 1 $, $ q \in \lbrace 1, 2 \rbrace $ and for any $ \beta \in \N^d $ with $ 0 \leq \abs{\beta} \leq 2 $,
	 \begin{align} \label{phi_n_bound}
	 	\left\| \partial_\beta \phi_n \right\|_q \leq \Cst{phi_n},
	 \end{align}
	 (note that, when $ \beta = (0, \ldots, 0) $, $ \partial_\beta \phi = \phi $).
	 
	 \subsection{The SLFV with mutations} \label{subsec:def_slfv}
	 
	 Let us now define the SLFV with mutations.
	 Fix $ u \in (0, 1] $, $ \mu > 0 $ and suppose that $ \nu(dr) $ is a finite measure on $ (0, \infty) $ satisfying
	 \begin{align} \label{condition_nu}
	 	\int_{0}^{\infty} r^d \nu(dr) < \infty.
	 \end{align}
	 Then, starting from an initial state $ \rho_0 \in \Xi $, the SLFV is defined as follows.
	 
	 \begin{definition}[The SLFV with mutations] \label{def:slfv}
	 	Let $ \Pi $ be a Poisson random measure on $ \R_+ \times \R^d \times (0,\infty) $ with intensity measure $ dt \otimes dx \otimes \nu(dr) $.
	 	For each point $ (t, x, r) \in \Pi $, a reproduction event takes place in the ball of centre $ x $ and radius $ r $ at time $ t $.
	 	At each reproduction event, we do the following:
	 	\begin{enumerate}
	 		\item choose a location $ y $ uniformly in $ B(x,r) $ and sample a parental type $ k_0 \in [0,1] $ according to the probability distribution $ \rho_{t^-}(y, dk) $,
	 		\item update $ \rho $ inside $ B(x,r) $ as follows:
	 		\begin{align} \label{update_rho}
	 			\forall z \in B(x,r), \quad \rho_t(z, dk) = (1-u) \rho_{t^-}(z, dk) + u \delta_{k_0}(dk).
	 		\end{align}
	 	\end{enumerate}
 	    Furthermore, for all $ s \leq t $ and $ x \in \R^d $ such that $ x $ does not find itself in the region affected by a reproduction event between times $ s $ and $ t $,
 	    \begin{align} \label{def_mutations}
 	    \rho_t(x,dk) = e^{-\mu (t-s)} \rho_s(x,dk) + \left( 1-e^{-\mu (t-s)} \right) dk,
 	    \end{align}
 		where $ dk $ denotes the Lebesgue measure on $ [0,1] $.
	 \end{definition}
 
 	 In other words, at each reproduction event, a proportion $ u $ of the individuals present in the ball $ B(x,r) $ dies and is replaced by the offspring of an individual sampled uniformly from inside this ball, while each individual, at rate $ \mu $, mutates to a new type sampled uniformly from the interval $ [0,1] $.
 	 Indeed, equation \eqref{def_mutations} says that, informally, between reproduction events, $ \rho_t $ solves the following
 	 \begin{align*} 
 	 \partial_t \rho_t(x, dk) = \mu(dk - \rho_t(x, dk)).
 	 \end{align*}
 	 The parameter $ u $ is called the \textit{impact parameter}, and $ \mu $ is called the mutation rate.
 	 
 	 The following remark will be crucial for the derivation of the Wright-Malécot formula.
 	 Note that, even though the SLFV is a continuous model in which there are no discrete units which correspond to ``individuals'', this model has an intrinsic genealogical structure, in the sense that the ancestry of random sample (of types) from the population can be described by a set of \emph{lineages}.
 	 These lineages form what is called a ``dual process'', and take the form of a system of coalescing random walks on $ \R^d $, see \cite{veber_spatial_2015}.
 	 When it finds itself in the region affected by a reproduction event, each lineage then jumps with probability $ u $ to a new location sampled uniformly in this region, and if several lineages jump during the same reproduction, they coalesce.
 	 
 	 This dual corresponds to the intuitive notion of ancestry which is implicit in Definition~\ref{def:slfv}.
 	 For the model of Definition~\ref{def:slfv}, a type sampled from the population at some location $ x $ is determined by sampling a lineage starting at $ x $, which is affected by mutations at rate $ \mu $.
 	 The sampled type is then determined either by the type sampled at the last mutation along the lineage, or by the distribution of ancestral types at the location of the lineage at time zero if no mutation has affected the lineage.
 	 Since the dual process of the SLFV is not used in the present work, we refer to \cite{barton_new_2010} and \cite{veber_spatial_2015} for a more rigorous discussion.
 	 
 	 \begin{remark} \label{remark:identity_in_state}
 	 	If we take $ \rho_0(x, dk) = dk $ for all $ x \in \R^d $, then, almost surely, at any time, two types sampled from the population are equal if and only if the corresponding lineages share a common ancestor at some point in the past and if neither of them has undergone a mutation since then.
 	 	The process $ (\rho_t, t \geq 0) $ can then be seen as tracking the size and geographical spread of all the ``families'' in the population, where a family is a macroscopic fraction of the population sharing the same type, \textit{i.e.} a portion of the population of the form $ f(x)\delta_{k_0}(dk) $ with $ f \geq 0 $ such that
 	 	\begin{align*}
 	 		\langle \rho_t - f\delta_{k_0}, \1_{k=k_0} \rangle = 0.
 	 	\end{align*}
 	 	One may see this as a generalisation of what is called \textit{tracer dynamics}, as introduced in \cite{hallatschek_gene_2008} (see also \cite{durrett_genealogies_2016}).
 	 \end{remark}
  
	 \begin{proposition} \label{prop:existence_uniqueness}
	 	There exists a unique $ \Xi $-valued Hunt process $ (\rho_t, t \geq 0) $ satisfying Definition~\ref{def:slfv}.
	 \end{proposition}
 	 
 	 This proposition follows directly from Corollary~2.4 in \cite{veber_spatial_2015} where the authors use a genealogical construction of the SLFV with mutations (see also \citep[Theorem~4.1]{etheridge_genealogical_2019}).
	 
	 \section{Main results} \label{sec:results}
	 
	 Let us now present our main results.
	 First, we obtain a central limit theorem for the SLFV with mutations of Definition~\ref{def:slfv} in two different regimes of reproduction events, each corresponding to a type of intensity measure $ \nu(dr) $.
	 The first case, called the ``fixed radius case'', corresponds to $ \nu(dr) $ being a Dirac measure at some fixed value $ R $.
	 In the second case, called the ``stable case'', we choose a measure $ \nu(dr) $ with a density which decays like a power of $ r $ as $ r \to\infty $.
	 In both cases, we rescale the SLFV in such a way that the measure $ \rho_t $ is very close to the Lebesgue measure on $ \R^d \times [0,1] $, which we denote by $ \lambda \in \Xi $, \textit{i.e.}
	 \begin{align} \label{def:lambda}
	 	\lambda(x, dk) = dk, \quad \forall x \in \R^d.
	 \end{align}
	 Furthermore, we rescale the difference $ \rho_t - \lambda $ so that it converges to a limiting process which we characterise as the solution to a stochastic partial differential equation (SPDE).
	 
	 From these results it then becomes possible to compute the asymptotic behaviour of the so-called probability of identity by descent, that is the probability that two individuals sampled from two prescribed locations share a common ancestor which is more recent than the last mutation to have affected their lineages.
	 
	 \subsection{The central limit theorem for the SLFV with mutations} \label{subsec:result_clt}
	 
	 We are interested in the asymptotic fluctuations of the SLFV with mutations from Definition~\ref{def:slfv} when both the impact parameter $ u $ and the mutation rate $ \mu $ tend to zero.
	 As we let these two parameters tend to zero, we also rescale space and time, so that the process is observed at times of the order of $ 1/\mu $ and over spatial scales of the order of $ \sqrt{u/\mu} $ (or a different power of this quantity in the stable case), assuming that $ u/\mu \to \infty $.
	 We state the result separately for the fixed radius case and the stable case, since the exact scalings differ slightly, but we shall provide a single proof covering both results in Section~\ref{sec:proof_clt}.
	 
	 Let us introduce some basic notations which will be used throughout the paper.
	 Let $ \mathcal{S}(\R^d \times [0,1]) $ denote the Schwartz space of rapidly decreasing smooth functions on $ \R^d \times [0,1] $, whose derivatives of all order are also rapidly decreasing.
	 More precisely, $ \phi \in \mathcal{S}(\R^d \times [0,1]) $ if, for any $ p \geq 1 $ and $ \beta \in \N^d $,
	 \begin{align*}
	 \sup_{(x, k) \in \R^d \times [0,1]} (1+\|x\|^p) \abs{ \partial_\beta \phi(x,k)} < \infty.
	 \end{align*}
	 Also if $ \phi $ and $ \psi $ are two functions defined on $ \R^d \times [0,1] $, we set
	 \begin{align*}
	 \phi \otimes \psi (x_1, k_1, x_2, k_2) = \phi(x_1, k_1) \psi(x_2, k_2).
	 \end{align*}
	 Accordingly, let $ \mathcal{S}'(\R^d \times [0,1]) $ denote the space of tempered distributions and let $ D(\R_+$, $\mathcal{S}'(\R^d \times [0,1])) $ denote the Skorokhod space of \cadlag distribution-valued processes (see Chapter~4 in \cite{walsh_introduction_1986}).
	 
	 \subsubsection{The fixed radius case} \label{subsubsec:clt_fixed_radius}
	 
	 Fix $ u \in (0,1] $, $ \mu > 0 $ and $ R > 0 $.
	 Let $ (\delta_N, N \geq 1) $ be a sequence of positive real numbers decreasing to zero and set, for $ N \geq 1 $,
	 \begin{align} \label{parameters_fixed_radius}
	 u_N = \frac{u}{N}, && \mu_N = \delta_N^2 \frac{\mu}{N}.
	 \end{align}
	 Further, for $ N\geq 1 $, let $ (\rho^N_t, t \geq 0) $ be the SLFV of Definition~\ref{def:slfv} with impact parameter $ u_N $, mutation rate $ \mu_N $ and with $ \nu(dr) = \delta_R(dr) $, started from $ \rho^N_0 = \lambda $.
	 Define the rescaled process $ (\bm{\rho}^N_t, t \geq 0) $ by setting, for $ N\geq 1 $,
	 \begin{align*}
	 \bm{\rho}^N_t(x, dk) := \rho^N_{Nt / \delta_N^2}\left(x / \delta_N, dk \right).
	 \end{align*}
	 Let $ \Delta $ denote the Laplace operator acting on the space variable, \textit{i.e.}, for $ \phi : \R^d \times [0,1] \to \R $ twice continuously differentiable in the space variable,
	 \begin{align*}
	 \Delta \phi(x, k) = \sum_{i=1}^{d} \deriv*[2]{\phi}{x_i}(x,k).
	 \end{align*}
	 Finally for $ r > 0 $, let $ V_r $ denote the volume of the $ d $-dimensional ball of radius $ r $.
	 The first important result of this paper is the following.
	 
	 \begin{theorem}[Central limit theorem for the SLFV with mutations - the fixed radius case] \label{thm:clt_fixed_radius}
	 	Suppose that $ \delta_N \to 0 $ and $ N \delta_N^{2-d} \to \infty $ as $ N \to \infty $ (note that the second condition is automatically satisfied when $ d \geq 2 $).
	 	Then, for all $ T > 0 $,
	 	\begin{align*}
	 	\lim_{N \to \infty} \E{ \sup_{t \in [0,T]} d(\bm{\rho}^N_t, \lambda) } = 0.
	 	\end{align*}
	 	Furthermore,
	 	\begin{align*}
	 	Z^N_t = (N \delta_N^{2-d})^{1/2} \left( \bm{\rho}^N_t - \lambda \right)
	 	\end{align*}
	 	defines a sequence of distribution-valued processes which converges in distribution in $ D(\R_+, \mathcal{S}'(\R^d \times [0,1])) $ to a process $ (Z_t, t \geq 0) $ which is the unique solution of the following SPDE:
	 	\begin{equation} \label{spde_Z_fixed_radius}
	 	\left\lbrace
	 	\begin{aligned}
	 	& d Z_t = \left[ \frac{\sigma^2}{2} \Delta Z_t - \mu Z_t \right] dt + d W(t), \\
	 	& Z_0 = 0,
	 	\end{aligned}
	 	\right.
	 	\end{equation}
	 	where $ \sigma^2 = u V_R \frac{2 R^2}{d+2} $ and $ (W(t), t\geq 0) $ is a Wiener process on $ \R^d \times [0,1] $ with covariation measure given by
	 	\begin{align} \label{def_Q}
	 	\mathcal{Q}(dx_1 dk_1 dx_2 dk_2) = u^2 V_R^2\, d x_1 \delta_{x_1}(dx_2) \left( dk_1 \delta_{k_1}(d k_2) - dk_1 dk_2 \right).
	 	\end{align}
	 \end{theorem}
	 
	 In other words, as the impact parameter and the mutation rate tend to zero according to \eqref{parameters_fixed_radius}, the rescaled SLFV converges to the uniform measure $ \lambda $ and the asymptotic deviations from this uniform measure are given by the process $ (Z_t, t \geq 0) $, where, for all $ \phi \in \mathcal{S}(\R^d \times [0,1]) $,
	 \begin{align*}
	 t \mapsto \langle Z_t, \phi \rangle - \int_{0}^{t} \left\langle Z_s, \frac{\sigma^2}{2} \Delta \phi - \mu \phi \right\rangle ds
	 \end{align*}
	 is a continuous square-integrable martingale with quadratic variation
	 \begin{align*}
	 t\, \langle \mathcal{Q}, \phi \otimes \phi \rangle = t\, u^2 V_R^2\, \int_{\R^d \times [0,1]} \left( \phi(x,k) - \int_{[0,1]} \phi(x,k')dk' \right)^2 dx dk.
	 \end{align*}
	 Note that the existence of $ W $ is stated in Proposition~\ref{prop:wiener_fixed_r}, and the above martingale problem uniquely characterises the distribution of the process $ (Z_t, t \geq 0) $, from \citep[Theorem~5.1]{walsh_introduction_1986}.
	 
	 Theorem~\ref{thm:clt_fixed_radius} is proved in Section~\ref{sec:proof_clt}.
	 The proof relies on a semimartingale decomposition of $ (\langle Z^N_t, \phi \rangle, t \geq 0) $, and a convergence theorem for sequences of stochastic integrals with respect to martingale measures.
	 
	 \begin{remark} \label{remark:generalisation_clt_fixed_radius}
	 	\begin{enumerate}
	 		\item Theorem~\ref{thm:clt_fixed_radius} is also true if we replace the measure $ \nu(dr) $ by any finite measure on $ (0,\infty) $ with a compact support.
	 		The coefficient $ \sigma^2 $ in \ref{spde_Z_fixed_radius} should then be replaced by
	 		\begin{align*}
	 			\sigma^2 = \frac{2 u}{d+2} \int_{0}^{\infty} r^2 V_r\, \nu(dr),
	 		\end{align*}
	 		and the coefficient $ u^2 V_R^2 $ in the covariation measure $ \mathcal{Q} $ should be replaced by
	 		\begin{align*}
	 			u^2 \int_{0}^{\infty} V_r^2\, \nu(dr).
	 		\end{align*}
	 		
	 		\item Different mutation mechanisms can also be considered. For example one could assume that, at each reproduction event, some proportion (say $ \mu $) of the offspring chooses a type uniformly in $ [0,1] $, or that all the offspring chooses a different type with some probability $ \mu $.
	 		We could also assume that mutants pick a new type according to a more general probability measure on $ [0,1] $, say $ \pi $.
	 		We would then need to replace $ \lambda $ by $ x \mapsto \pi(dk) $ in the definition of $ Z^N $, and the covariation measure $ \mathcal{Q} $ should be replaced by
	 		\begin{align*}
	 			u^2 V_R^2\, dx_1 \, \delta_{x_1}(dx_2) \, (\pi(dk_1)\, \delta_{k_1}(dk_2) - \pi(dk_1)\, \pi(dk_2)).
	 		\end{align*}
	 	\end{enumerate}
	 \end{remark}
 
 	Since the main purpose of this paper is to study the stationary behaviour of the SLFV with mutations, we have assumed that $ \bm{\rho}^N_0 = \lambda $ for all $ N \geq 1 $.
 	In Section~\ref{subsec:non-stationary}, we state a more general version of Theorem~\ref{thm:clt_fixed_radius} (and of Theorem~\ref{thm:clt_stable} below) where the initial condition $ \bm{\rho}^N_0 $ is assumed to converge to some $ \bm{\rho}_0 \in \Xi $.
 	In addition, we also consider a more general mutation mechanism, where the trait distribution at each spatial location evolves according to a Feller semigroup on $ [0,1] $ between reproduction events (which amounts to assuming that the trait of each individual evolves according to a Markov process on $ [0,1] $ during its lifetime).

     \subsubsection{The stable case} \label{subsubsec:clt_stable}
     
     We now want to extend the previous analysis to a situation in which reproduction events can affect arbitrarily large regions and such that these large scale reproduction events take place often enough to significantly alter the qualitative behaviour of the SLFV.
     This will result in increased correlations between the genetic compositions of different spatial locations, both through non-local diffusion and correlations in the noise driving the fluctuations of the limiting process.
     
     First fix $ \alpha \in (0, d \wedge 2) $ and set
     \begin{align*}
     \nu_\alpha(dr) := \frac{\1_{r \geq 1}}{r^{d+\alpha +1}} dr.
     \end{align*}
     It is straightforward to check that $ \nu_\alpha $ satisfies \eqref{condition_nu}.
     Also fix $ u \in (0,1] $ and $ \mu > 0 $ and let $ (\delta_N, N \geq 1) $ be a sequence of positive numbers decreasing to zero.
     For $ N \geq 1 $, set
     \begin{align} \label{scaling_stable}
     u_N = \frac{u}{N}, && \mu_N = \delta_N^\alpha \frac{\mu}{N},
     \end{align}
     and let $ (\rho^N_t, t \geq 0) $ be the SLFV of Definition~\ref{def:slfv} with impact parameter $ u_N $, mutation rate $ \mu_N $ and with $ \nu = \nu_\alpha $, started from $ \rho_0 = \lambda $.
     Define the rescaled SLFV as
     \begin{align*}
     \bm{\rho}^N_t(x, dk) = \rho^N_{Nt / \delta_N^\alpha}\left( x / \delta_N, dk \right).
     \end{align*}
     
     Before stating our result, we introduce some notations.
     First, for $ x, y \in \R^d $, set
     \begin{align} \label{def_V_2r}
     V_{2,r}(x,y) = \int_{\R^d} \1_{\lbrace \|x-z\| < r, \|y - z\| < r \rbrace} dz
     \end{align}
     and
     \begin{align*}
     \Phi_\alpha(\|x-y\|) = \int_{\frac{\|x-y\|}{2}}^{\infty} \frac{V_{2,r}(x,y)}{V_r} \frac{dr}{r^{d+\alpha+1}}.
     \end{align*}
     Let $ \mathcal{D}_\alpha $ be the fractional Laplacian acting on the space variable, \textit{i.e.} for any $ \phi : \R^d \times [0,1] \to \R $ admitting uniformly bounded spatial derivatives of order at least two (see \cite{samko_fractional_1993}),
     \begin{align} \label{def_Dalpha}
     \mathcal{D}_\alpha \phi(x,k) = \int_{\R^d} \Phi_\alpha(\|x-y\|) (\phi(y,k) - \phi(x,k) - \1_{\lbrace \|x-y\| \leq 1 \rbrace} \nabla \phi(x,k) \cdot (y-x)) dy.
     \end{align}
     For $ x, y \in \R^d $, also set
     \begin{align} \label{def_Kalpha}
     K_\alpha(x, y) = \int_{\frac{\|x-y\|}{2}}^{\infty} V_{2,r}(x,y) \frac{dr}{r^{d+\alpha + 1}} = \frac{C_{d,\alpha}}{\|x-y\|^\alpha}
     \end{align}
     where $ C_{d,\alpha} $ is a positive constant depending only on $ d $ and $ \alpha $.
     
     We now state our second main result.
     
     \begin{theorem}[Central limit theorem for the SLFV with mutations - the stable case] \label{thm:clt_stable}
     	Assume that $ \delta_N \to 0 $ as $ N \to \infty $.
     	Then, for all $ T > 0 $,
     	\begin{align*}
     	\lim_{N \to \infty} \E{ \sup_{t \in [0,T]} d(\bm{\rho}^N_t, \lambda) } = 0.
     	\end{align*}
     	Furthermore,
     	\begin{align*}
     	Z^N_t = \sqrt{N} ( \bm{\rho}^N_t - \lambda )
     	\end{align*}
     	defines a sequence of distribution-valued processes which converges in distribution in $ D(\R_+,$ $\mathcal{S}'(\R^d \times [0,1])) $ to a process $ (Z_t, t \geq 0) $ which is the unique solution to the following SPDE:
     	\begin{equation} \label{spde_Z_stable}
     	\left\lbrace
     	\begin{aligned}
     	& d Z_t = \left[ u \mathcal{D}_\alpha Z_t - \mu Z_t \right] dt + d W(t) \\
     	& Z_0 = 0,
     	\end{aligned}
     	\right.
     	\end{equation}
     	where $ (W(t), t \geq 0) $ is a Wiener process on $ \R^d \times [0,1] $ with covariation measure on $ (\R^d \times [0,1])^2 $ given by
     	\begin{align} \label{def_Qalpha}
     	\mathcal{Q}_\alpha (dx_1 dk_1 dx_2 dk_2) = u^2\, K_\alpha(x_1, x_2) dx_1 dx_2 ( dk_1 \delta_{k_1}(dk_2) - dk_1 dk_2).
     	\end{align}
     \end{theorem}
     
     The martingale problem associated to \eqref{spde_Z_stable} is the following.
     For any $ \phi \in \mathcal{S}(\R^d \times [0,1]) $,
     \begin{align*}
     t \mapsto \langle Z_t, \phi \rangle - \int_{0}^{t} \left\langle Z_s, u\mathcal{D}_\alpha \phi - \mu \phi \right\rangle ds
     \end{align*}
     is a square integrable continuous martingale with quadratic variation
     \begin{multline*}
     t\, \langle \mathcal{Q}_\alpha, \phi \otimes \phi \rangle = t\, u^2 \int_{(\R^d)^2 \times [0,1]} \left( \phi(x_1, k) - \int_{[0,1]} \phi(x_1,k')dk' \right) \\ \times \left( \phi(x_2, k) - \int_{[0,1]} \phi(x_2, k')dk' \right) \frac{C_{d,\alpha}}{\| x_1 - x_2 \|^\alpha} dx_1 dx_2 dk.
     \end{multline*}
     We recall why there exists such a $ W $ in Proposition~\ref{prop:wiener_stable}, and as before, there exists a unique process $ (Z_t, t \geq 0) $ solving \eqref{spde_Z_stable}, by \citep[Theorem~5.1]{walsh_introduction_1986}.
     
     Theorem~\ref{thm:clt_stable} is proved along with Theorem~\ref{thm:clt_fixed_radius} in Section~\ref{sec:proof_clt}.
     The main differences with Theorem~\ref{thm:clt_fixed_radius} are that the Laplacian is replaced with the non-local operator $ \mathcal{D}_\alpha $ and that the Gaussian noise driving the fluctuations is now correlated in space, with correlations decaying as $ \|x-y\|^{-\alpha} $.
     These two changes result from the large scale reproduction events which take place rarely enough that the population retains a signature of isolation by distance but often enough to induce these strong spatial correlations.
     The index $ \alpha $ is a convenient measure of the strength of these correlations: the closer it is to zero the stronger they are and the closer it is to 2 the more localised the correlations become.
	 
	 \subsection{The Wright-Malécot formula for isolation by distance under short range and long range dispersal} \label{subsec:result_wright_malecot}
	 
	 Despite their apparent complexity and lack of direct reference to ancestry, Theorems~\ref{thm:clt_fixed_radius} and~\ref{thm:clt_stable} are deeply linked to previous results on the sharing of recent common ancestors in a spatially distributed population.
	 These include results on the stepping stone model \cite{kimura_stepping_1964,sawyer_asymptotic_1977} and the SLFV \cite{barton_inference_2013}, also see \cite{barton_neutral_2002}.
	 Indeed, Theorems~\ref{thm:clt_fixed_radius} and~\ref{thm:clt_stable} can be seen as results on the correlations between the genetic compositions of the population at different spatial locations.
	 
	 To see this, consider the following.
	 Let $ \phi $ and $ \psi $ be two probability density functions on $ \R^d $.
	 Sample two locations $ x_1 $, $ x_2 $ according to $ \phi $ and $ \psi $, respectively, and sample one genetic type at each of these locations according to the distribution of types in $ \bm{\rho}^N_t $ at some time $ t \geq 0 $.
	 Let $ P^N_t(\phi, \psi) $ be the probability that these two types are the same.
	 
	 In the vocabulary of population genetics, $ P^N_t(\phi, \psi) $ is the probability of identity in state of two ``individuals'' sampled according to $ \phi $ and $ \psi $.
	 In view of Remark~\ref{remark:identity_in_state}, this coincides with the probability of identity by descent, \textit{i.e.} the probability that the two sampled lineages share a common ancestor that is more recent than the last time either lineage experienced a mutation.
	 
	 This probability can be written more explicitly in terms of the process $ (\bm{\rho}^N_t, t \geq 0) $ as follows.
	 Let $ \1_{\Delta} : [0,1]^2 \to \R $ denote the indicator function of the diagonal, \textit{i.e.}
	 \begin{align*}
	 	\1_{\Delta}(k_1, k_2) = \1_{k_1 = k_2}.
	 \end{align*}
	 Then,
	 \begin{align} \label{definition_PN}
	 	P^N_t(\phi, \psi) = \E{ \left\langle \bm{\rho}^N_t \otimes \bm{\rho}^N_t, (\phi \otimes \psi) \1_{\Delta} \right\rangle },
	 \end{align}
	 where
	 \begin{align*}
	 	(\phi \otimes \psi) \1_{\Delta} (x_1, k_1, x_2, k_2) = \phi(x_1) \psi(x_2) \1_{k_1 = k_2}.
	 \end{align*}
	 
	 The following is then a consequence of Theorems~\ref{thm:clt_fixed_radius} and~\ref{thm:clt_stable}.
	 Let $ G^{(\alpha)}_t : \R^d \to \R $ denote the fundamental solution associated to $ \partial_t - \mathcal{D}_\alpha $, \textit{i.e.} such that
	 \begin{align} \label{def_G_alpha}
	 \partial_t G_t^{(\alpha)} = \mathcal{D}_\alpha G^{(\alpha)}_t
	 \end{align}
	 and
	 \begin{align*}
	 \int_{\R^d} G^{(\alpha)}_t (x-y) \phi (y) dy \cvgas{t}{0} \phi(x),
	 \end{align*}
	 for any twice continuously differentiable $ \phi : \R^d \to \R $.
	 
	 \begin{theorem}[Wright-Malécot formula for identity by descent] \label{thm:WM_formula}
	 	Assume that $ \phi $ and $ \psi $ are two smooth and compactly supported probability density functions on $ \R^d $. 
	 	Then, under the conditions of Theorem~\ref{thm:clt_fixed_radius},
	 	\begin{align} \label{Wright-Malecot_fixed_radius}
	 		\lim_{t \to \infty} \lim_{N \to \infty} N \delta_N^{2-d} P^N_t(\phi, \psi) = \frac{u^2 V_R^2}{(2\pi \sigma^2)^{d/2}} \int_{(\R^d)^2} F \left( \frac{\|x-y\|}{\sigma} \right) \phi(x) \psi(y) dx dy,
	 	\end{align}
	 	where $ \sigma^2 = u V_R \frac{2 R^2}{d + 2} $ and the function $ F $ depends only on $ d $ and $ \mu $ and is given by
	 	\begin{align*}
	 		F(x) = \left( \frac{x}{\sqrt{2\mu}} \right)^{1 - d / 2} K_{1 - d/2} \left( \sqrt{2\mu} x \right),
	 	\end{align*}
	 	where $ K_\nu $ denotes the modified Bessel function of the second kind of degree $ \nu $ \cite{abramowitz_handbook_1964}.
	 	On the other hand, under the conditions of Theorem~\ref{thm:clt_stable},
	 	\begin{align} \label{Wright-Malecot_stable}
	 	\lim_{t \to \infty} \lim_{N \to \infty} N P^N_t(\phi, \psi) = u \int_{(\R^d)^2} F_{d, \alpha} \left( \left(\mu / u\right)^{1/\alpha} \|x-y\| \right) \phi(x) \psi(y) dx dy,
	 	\end{align}
	 	where $ F_{d, \alpha} : \R_+ \to \R_+ $ is such that, for any $ x, y \in \R^d $,
	 	\begin{align*}
	 	F_{d,\alpha}(\|x-y\|) = \int_{0}^{\infty} \int_{(\R^d)^2} e^{-2t} G^{(\alpha)}_t(x-z_1) G^{(\alpha)}_t(y-z_2) \frac{C_{d,\alpha}}{\|z_1 - z_2\|^\alpha} dz_1 dz_2 dt.
	 	\end{align*}
	 	(The fact that the right hand side only depends on $ \|x-y\| $ can be seen by a change of variables.)
	 \end{theorem}
 
\begin{figure}[h]
 	 	\centering
 	 	\includegraphics[height=0.35\textheight]{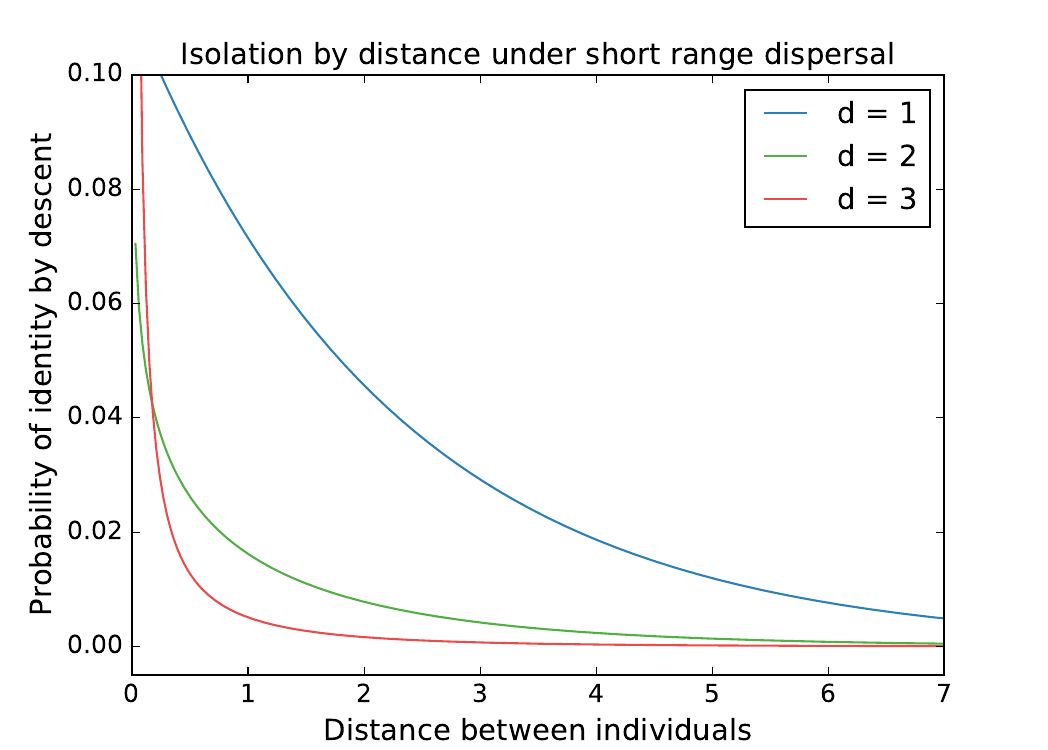}\par
 	 	\caption{Deacrease of the probability of identity by descent as a function of the distance between individuals. The graph shows a plot of the function $ F $ for $ d \in \lbrace 1, 2, 3 \rbrace $. We can see that it decreases much more rapidly for larger values of $ d $ than for smaller values.}
 	 	\label{fig:wright_malecot_fixed_radius}
 	 \end{figure}

 We prove Theorem~\ref{thm:WM_formula} in Section~\ref{sec:WM_proof}.
 	 The convergence \eqref{Wright-Malecot_fixed_radius} should be compared to equations (10) and (15) in \cite{barton_neutral_2002} (originally due to Malécot \cite{malecot_heterozygosity_1975}) or (1.13) and (2.22) in \cite{kimura_stepping_1964}.
 	 This is known in the literature as the Wright-Malécot formula, and is widely used to infer both the mean-square displacement of individuals in the population (\textit{i.e.} $ \sigma^2 $) and the effective population density from genetic samples \cite{rousset_genetic_1997}.
 	 Figure~\ref{fig:wright_malecot_fixed_radius} shows the behaviour of the function $ F $ for $ d \in \lbrace 1, 2, 3 \rbrace $.
 	 
 	 \begin{remark} \label{remark:degeneracy}
 	 	It can be noted that for $ d \geq 2 $, the function $ F $ is degenerate at zero, \textit{i.e.}
 	 	\begin{align*}
 	 		\lim_{x \downarrow 0} F(x) = + \infty.
 	 	\end{align*}
 	 	This can be surprising since $ P^N_t(\phi, \psi) $ was defined as a probability.
 	 	This reflects the fact that the Wright-Malécot formula breaks down if we try to sample two individuals from \textit{exactly} the same location.
 	 	That is why we need to integrate against the probability density functions $ \phi $ and $ \psi $ (it can be shown from the proof of Theorem~\ref{thm:WM_formula} that the right hand side of \eqref{Wright-Malecot_fixed_radius} is bounded by a constant times $ \| \phi \|_2 \| \psi \|_2 $).
 	 \end{remark}
 	 
 	 The second part of Theorem~\ref{thm:WM_formula} is a very important step towards developing statistical inference procedures from genetic data adapted to species undergoing long range dispersal.
 	 Indeed, with this result it could become possible to estimate the parameters $ u $ and $ \alpha $ from a sample of genetic markers from different spatial locations in the population.
 	 One hurdle which remains in the way is that we have to find an efficient way to compute numerically the function $ F_{d,\alpha} $.
 	 This might be done using Fourier transforms, but the heavy tail of the function $ G^{(\alpha)}_t(\cdot) $ makes any rigorous analysis quite daunting.
   	 
   	 \begin{remark} \label{remark:universality_WM}
   	 	 The Wright-Malécot formula has been shown to hold for a wide variety of spatial models in population genetics \cite{barton_neutral_2002,barton_inference_2013}.
   	 	 The proof of Theorem~\ref{thm:WM_formula} shows that this formula is directly linked to the limiting behaviour of the fluctuations in the genetic composition of the population, and that any model which displays the same asymptotical behaviour as the SLFV in Theorem~\ref{thm:clt_fixed_radius} should satisfy the Wright-Malécot formula.
   	 \end{remark}
    
    \subsection{The SLFV with general mutation mechanism and in the non-stationary regime} \label{subsec:non-stationary}
    
    Theorems~\ref{thm:clt_fixed_radius} and \ref{thm:clt_stable} allowed us to derive the Wright-Malécot formula from the Gaussian fluctuations of the SLFV of Definition~\ref{def:slfv} (Theorem~\ref{thm:WM_formula}).
    The aim of this section is to restate Theorems~\ref{thm:clt_fixed_radius} and \ref{thm:clt_stable} with more general assumptions on the mutation mechanism and on the initial condition of the process.
    For simplicity, we still assume that mutations take place during an individual's lifetime, and not at birth.
    We shall suppose that the type of an individual follows a Markov process on $ [0,1] $.
    Since there are no individuals \textit{per se} in the model, this translates to the following definition.
    
    \begin{definition}[The SLFV with general mutation mechanism] \label{def:slfv_gl_mutations}
    	Suppose that $ (\mathcal{T}_s, s \geq 0) $ is a Feller semigroup acting on bounded measurable real-valued functions defined on $ [0,1] $ (hence $ \mathcal{T}^*_s $ acts on $ \mathcal{M}_1([0,1]) $).
    	We then define a process $ (\rho_t, t \geq 0) $ which evolves exactly as in Definition~\ref{def:slfv}, except that, for all $ s \leq t $ and $ x \in \R^d $ such that $ x $ does not find itself in the region affected by a reproduction event between times $ s $ and $ t $,
    	\begin{align*}
    		\rho_t(x,\cdot) = \mathcal{T}_{\mu (t-s)}^* \rho_s(x, \cdot).
    	\end{align*}
    \end{definition}
    
    Let $ \mathcal{G} $ denote the generator of the semigroup $ (\mathcal{T}_s, s \geq 0) $, whose domain will be denoted by $ \mathcal{D}(\mathcal{G}) $.
    We shall consider the same scaling as before, \textit{i.e.} we assume that
    \begin{align*}
    	u_N = \frac{u}{N}, && \mu_N = \delta_N^\alpha \frac{\mu}{N},
    \end{align*}
    where $ \alpha = 2 $ in the fixed radius case, and we suppose that $ (\rho^N_t, t \geq 0) $ is the SLFV of Definition~\ref{def:slfv_gl_mutations} with impact parameter $ u_N $, mutation rate $ \mu_N $, and with $ \nu = \delta_R $ in the fixed radius case, and $ \nu = \nu_\alpha $ in the stable case.
    Then let
    \begin{align*}
    	\bm{\rho}^N_t(x,dk) := \rho^N_{Nt / \delta_N^\alpha} (x / \delta_N, dk).
    \end{align*}
    Recall the definition of $ G^{(\alpha)}_t $ in \eqref{def_G_alpha} and let $ G^{(2)}_t : \R^d \to \R $ be the fundamental solution associated to $ \partial_t - \frac{\sigma^2}{2u} \Delta $, \textit{i.e.}
    \begin{align} \label{def_G_2}
    	G^{(2)}_t(x) = \frac{1}{\left(2\pi t\, \sigma^2/u\right)^{d/2}} \exp \left( - \frac{\| x \|^2}{2 t\, \sigma^2/u} \right),
    \end{align}
    where $ \sigma^2 $ is as in Theorem~\ref{thm:clt_fixed_radius}.
    Assume that $ \bm{\rho}^N_0 $ converges weakly to a deterministic $ \bm{\rho}_0 \in \Xi $ in probability as $ N \to \infty $.
    We then define $ \bm{p}^{(\alpha)} = (\bm{p}_t^{(\alpha)}, t \geq 0) $ by
    \begin{align} \label{def_p_alpha}
    	\bm{p}^{(\alpha)}_t(x, dk) := \int_{\R^d} G^{(\alpha)}_{u t}(x-y) (\mathcal{T}^*_{\mu t}\, \bm{\rho}_0)(y, dk)\, dy.
    \end{align}
    The following generalises Theorems~\ref{thm:clt_fixed_radius} and \ref{thm:clt_stable} to the present setting.

    \begin{theorem} \label{thm:clt_non_stationary}
    	Let $ \eta_N = \delta_N^{2-d} $ in the fixed radius case and $ \eta_N = 1 $ in the stable case, and assume that $ N \eta_N \to \infty $ as $ N \to \infty $.
    	Also assume that $ \bm{\rho}^N_0 $ converges weakly (as a measure via the correspondence \eqref{xi_measure_oto_correspondence}) as $ N \to \infty $ to a deterministic $ \bm{\rho}_0 \in \Xi $ in probability.
    	Then, for all $ T > 0 $,
    	\begin{align} \label{cvg_det_limit_general}
    		\lim_{N \to \infty} \E{ \sup_{t \in [0,1]} d(\bm{\rho}^N_t, \bm{p}^{(\alpha)}_t) } = 0.
    	\end{align}
    	Furthermore, there exists a deterministic sequence $ (\bm{p}^N, N \geq 1) $ of elements of $ D(\R_+, \Xi) $, which converges (locally uniformly in time) to $ \bm{p}^{(\alpha)} $, such that
    	\begin{align*}
    		Z^N_t := (N \eta_N)^{1/2} (\bm{\rho}^N_t - \bm{p}^N_t)
    	\end{align*}
    	converges in distribution in $ D(\R_+, \mathcal{S}'(\R^d \times [0,1])) $ to a process $ (Z_t, t \geq 0) $.
    	In the fixed radius case, for any $ \phi \in \mathcal{S}(\R^d \times [0,1]) $ such that $ k \mapsto \phi(x,k) \in \mathcal{D}(\mathcal{G}) $ for all $ x\in \R^d $,
    	\begin{align*}
    		t \mapsto \langle Z_t, \phi \rangle - \int_{0}^{t} \left\langle Z_s, \frac{\sigma^2}{2} \Delta \phi + \mu \mathcal{G} \phi \right\rangle ds
    	\end{align*}
    	is a square-integrable continuous martingale with quadratic variation
    	\begin{align*}
    		(u V_R)^2 \int_{0}^{t} \int_{\R^d \times [0,1]} \left( \phi(x,k) - \int_{[0,1]} \phi(x,k') \bm{p}^{(2)}_s(x, dk') \right)^2 \bm{p}^{(2)}_s(x,dk) dx ds.
    	\end{align*}
    	In the stable case, on the other hand,
    	\begin{align*}
    		t \mapsto \langle Z_t, \phi \rangle - \int_{0}^{t} \left\langle Z_s, u \mathcal{D}_\alpha \phi + \mu \mathcal{G} \phi \right\rangle ds
    	\end{align*}
    	is a square-integrable continuous martingale with quadratic variation
    	\begin{align} \label{qvar_non_stationary}
    		u^2 \int_{0}^{t} \int_{0}^{\infty} \langle \bm{p}^{(\alpha)}_s, \Theta_{r,s} \rangle \frac{dr}{r^{1+d+\alpha}} ds,
    	\end{align}
    	where
    	\begin{align*}
    		\Theta_{r,s}(x,k) := \frac{1}{V_r} \int_{B(x,r)} \left[ \int_{B(y,r)} \left( \phi(z,k) - \int_{[0,1]} \phi(z,k') \bm{p}^{(\alpha)}_s(z, dk') \right) dz \right]^2 dy.
    	\end{align*}
    \end{theorem}

	The proof of Theorem~\ref{thm:clt_non_stationary} is almost identical to that of Theorems~\ref{thm:clt_fixed_radius} and \ref{thm:clt_stable}.
	Details on the specific adaptations which are needed in the proof are given in Section~\ref{subsec:proof_clt_non-stationary}.

    The main differences with Theorems~\ref{thm:clt_fixed_radius} and \ref{thm:clt_stable} are that the deterministic limit of $ \bm{\rho}^N $ and the centring term $ \bm{p}^N $ are no longer constant. 
    This affects the limiting fluctuations trough the covariation measure of the limiting driving noise, which now depends on $ \bm{p}^{(\alpha)} $.
    Furthermore, the term $ -\mu Z_t $ has now been replaced by $ \mu \mathcal{G}^* Z_t $ in the limiting SPDE, which accounts for the more general mutation mechanism.
    We note that, for the mutation mechanism considered in Theorems~\ref{thm:clt_fixed_radius} and \ref{thm:clt_stable},
    \begin{align*}
    	\mathcal{G} \phi (k) = \int_{[0,1]} (\phi(k')-\phi(k)) dk'.
    \end{align*}
    Hence $ \mathcal{G}^* Z_t \neq - Z_t $, but we see that, if $ \phi(x,k) = \psi(x) $, then $ \langle Z_t, \phi \rangle = 0 $ almost surely, and it follows that $ \langle \mathcal{G}^* Z_t, \phi \rangle = - \langle Z_t, \phi \rangle $, hence the two equations admit the same solution.
    
    Note that the deterministic limit $ \bm{p}^{(\alpha)} $ has the following interpretation: if $ (X^{(\alpha)}_t, t \geq 0) $ is a Markov process on $ \R^d $ with generator $ u D_\alpha $ in the stable case and $ \frac{\sigma^2}{2} \Delta $ in the fixed radius case, and if $ (\mathcal{K}_t, t \geq 0) $ is a Markov process on $ [0,1] $ with generator $ \mu \mathcal{G} $, then, for any smooth and compactly supported $ \phi : \R^d \times [0,1] \to \R $,
    \begin{align*}
    	\langle \bm{p}^{(\alpha)}_t, \phi \rangle = \int_{\R^d \times [0,1]} \E[(x,k)]{\phi(X^{(\alpha)}_t, \mathcal{K}_t)} \bm{\rho}_0(x, dk) dx,
    \end{align*}
    where $ \E[(x,k)]{\cdot} $ denotes the expectation with respect to the joint distribution of the pair of independent processes $ (X^{(\alpha)}, \mathcal{K}) $ started from $ (x,k) $.
    
    \begin{remark} \label{rk:centring_term}
    	In classical central limit theorems, the centring term is usually equal to the deterministic limit of the process whose fluctuations we are interested in.
    	Here we introduce a sequence of centring terms $ (\bm{p}^N, N \geq 1) $ because, while $ \bm{\rho}^N_t - \bm{p}^N_t $ is of the order of $ (N\eta_N)^{-1/2} $, the difference $ \bm{p}^N - \bm{p}^{(\alpha)} $ is at least of the order of $ (\delta_N)^{\gamma} $, where $ \gamma = 2 $ in the fixed radius case and $ \gamma = 2-\alpha $ in the stable case (see Proposition~\ref{prop:cvg_centring_term}).
    	It follows that the centring term in the definition of $ Z^N $ can be replaced by $ \bm{p}^{(\alpha)} $ only when $ (N \eta_N)^{1/2} (\delta_N)^\gamma \to 0 $ as $ N \to \infty $ and if $ \bm{\rho}^N_0 $ converges to $ \bm{\rho}_0 $ sufficiently quickly (and uniformly over $ \R^d $).
    	The condition on $ \delta_N $ translates to $ N (\delta_N)^{6-d} \to 0 $ in the fixed radius case and $ \sqrt{N} (\delta_N)^{2-\alpha} \to 0 $ in the stable case.
    	A similar issue arose in \cite{forien_central_2017}, where a similar sequence of centring terms had to be defined.
    \end{remark}
	 
	\section{Proof of the central limit theorems} \label{sec:proof_clt}
	 
	The proof of Theorems~\ref{thm:clt_fixed_radius}, \ref{thm:clt_stable} and \ref{thm:clt_non_stationary} is similar in spirit to what was done in \cite{forien_central_2017}: we write $ \langle \rho^N_t, \phi \rangle $ as the sum of a predictable term and a martingale term, and we use martingale convergence theorems to show tightness and convergence in distribution of $ (\langle Z^N_t, \phi \rangle, t \geq 0) $ for any $ \phi \in \mathcal{S}(\R^d \times [0,1]) $.
	It is then easy to generalise this to the joint convergence of $ (\langle Z^N_t, \phi_1 \rangle, \ldots, \langle Z^N_t, \phi_n \rangle) $ and obtain the convergence in distribution of $ (Z^N_t, t \geq 0) $ in $ D(\R_+, \mathcal{S}'(\R^d \times [0,1])) $.
	
	In the next subsection, we outline the strategy of the proof of the first two central limit theorems, stating a semimartingale decomposition of the process $ (Z^N_t, t \geq 0) $.
	This decomposition is proved in Subsection~\ref{subsec:semimartingale}.
	Subsections~\ref{subsec:bound_qvar} through \ref{subsec:convergence_Mn} then prepare several intermediary results needed for the proof.
	In particular, the convergence of $ (\bm{\rho}^N_t, t \in [0,T]) $ to its deterministic limit $ \lambda $ is proved in Subsection~\ref{subsec:det_limit}.
	The two central limit theorems are then proved in Subsection~\ref{subsec:proof_clt}, while the small adjustments needed for the proof of Theorem~\ref{thm:clt_non_stationary} are detailed in Subsection~\ref{subsec:proof_clt_non-stationary}.
	 
	 \subsection{Outline of the proof} \label{subsec:strategy}
	 
	 To avoid repetitions, we prove both Theorem~\ref{thm:clt_fixed_radius} and \ref{thm:clt_stable} at the same time, using general notations.
	 The parameter $ \alpha $ will thus be set equal to 2 in the fixed radius case, in agreement with \eqref{parameters_fixed_radius} and \eqref{scaling_stable}.
	 Recall the definition of $ E_q $ in Section~\ref{subsec:state_space}.
	 For $ \phi \in E_1 $ and $ r > 0 $, let $ \overline{\phi}(\cdot, r) $ be defined by
	 \begin{align} \label{def_average}
	 \overline{\phi}(x, k, r) := \frac{1}{V_r} \int_{B(x,r)} \phi(y, k) dy.
	 \end{align}
	 Then, $ \overline{\overline{\phi}}(x,k,r) $ denotes the average of the average of $ \phi $, \textit{i.e.}
	 \begin{align*}
	 	\overline{\overline{\phi}}(x,k,r) := \frac{1}{V_r^2} \int_{B(x,r)} \int_{B(y,r)} \phi(z,k) dz\, dy.
	 \end{align*}
	 We then set, for any $ \phi \in E_1 $,
	 \begin{align} \label{def_LN}
	 \mathcal{L}^{N, \alpha} \phi(x, k) := \delta_N^{-\alpha} \int_{0}^{\infty} V_r \left( \overline{\overline{\phi}}(x, k, \delta_N r) - \phi(x,k) \right) \nu_\alpha(dr),
	 \end{align}
	 where $ \nu_\alpha = \delta_R $ for $ \alpha = 2 $.
	 Note that the operator $ u\mathcal{L}^{N,\alpha} $ is the generator of a continuous-time random walk $ (X^{N,\alpha}_t, t \geq 0) $ on $ \R^d $ which jumps at rate $ u\, \delta_N^{-\alpha} \int_{0}^{\infty} V_r\, \nu_\alpha(dr) $ from its current position $ x $ to $ x + \delta_N \mathcal{R}(Y_1 + Y_2) $, where $ \mathcal{R} $ is a random variable with distribution
	 \begin{align*}
	 	\frac{V_r\, \nu_\alpha(dr)}{\int_{0}^{\infty} V_{r'}\, \nu_\alpha(dr')},
	 \end{align*}
	 and $ Y_1 $ and $ Y_2 $ are two independent uniform random variables on $ B(0,1) \subset \R^d $, also independent from $ \mathcal{R} $.
	 In the following, $ (P^{N,\alpha}_t, t \geq 0) $ denotes the strongly continuous semigroup generated by $ u \mathcal{L}^{N,\alpha}-\mu $, \textit{i.e.}
	 \begin{align} \label{def_PNalpha}
	 	P^{N,\alpha}_t\phi(x,k) := e^{-\mu t}\, \E[x]{\phi(X^{N,\alpha}_t, k)},
	 \end{align}
	 for any $ \phi : \R^d \times [0,1] \to \R $ bounded and measurable, where $ \E[x]{\cdot} $ denotes the expectation with respect to the distribution of $ (X^{N,\alpha}_t, t \geq 0) $, started from $ X^{N,\alpha}_0 = x $.
	 As we shall see below (Lemma~\ref{lemma:cvg_PN}), $ P^{N,\alpha}_t $ is a contraction from $ E_q $ into itself.
	 
	 For $ \rho \in \Xi $ and $ r > 0 $, we also define a map $ [\rho]_r : \R^d \times \R^d \to \mathcal{M}_1([0,1]) $ as follows, recalling the definition of $ V_{2,r}(x_1,x_2) $ in \eqref{def_V_2r},
	 \begin{equation} \label{def_bracket_rho}
	 [\rho]_r(x_1, x_2, dk) := \left\lbrace
	 \begin{aligned}
	 & \dfrac{1}{V_{2,r}(x_1, x_2)} \int_{B(x_1,r) \cap B(x_2,r)} \frac{1}{V_r} \int_{B(y,r)} \rho(z, dk) dz dy & \text{ if } \abs{x_1 - x_2} < 2r, \\
	 & 0 & \text{ otherwise.} \hspace{1.2cm}
	 \end{aligned}
	 \right.
	 \end{equation}
	 If $ \nu $ is a finite measure on $ (0,\infty) $ satisfying \eqref{condition_nu} and $ \rho \in \Xi $, we define a map $ \Gamma^\nu(\rho) : \R^d \times \R^d \to \mathcal{M}([0,1]^2) $ as follows
	 \begin{multline} \label{def_Gamma}
	 \Gamma^\nu(\rho)(x_1, x_2, dk_1 dk_2) := \int_{0}^{\infty} V_{2,r}(x_1, x_2) \Big[ [\rho]_r(x_1, x_2, dk_1) \delta_{k_1}(dk_2) - \rho(x_1, dk_1) [\rho]_r(x_1, x_2, dk_2)  \\ - [\rho]_r(x_1, x_2, dk_1) \rho(x_2, dk_2) + \rho(x_1, dk_1) \rho(x_2, dk_2) \Big] \nu(dr).
	 \end{multline}
	 When $ \nu = \delta_r $, we also write $ \Gamma^\nu = \Gamma^r $.
	 Finally, for $ N \geq 1 $, define the measure $ \nu^N_\alpha $ by
	 \begin{align} \label{def_nu_N_alpha}
	 \int_{0}^{\infty} f(r) \nu_\alpha^N(dr) = \int_{0}^{\infty} f(\delta_N r) \delta_N^{-(\alpha+d)} \nu_\alpha(dr).
	 \end{align}
	 The following proposition then gives the semimartingale form of the rescaled fluctuations process $ Z^N $.
	 It will be proved in Subsection~\ref{subsec:semimartingale}.
	 
	 \begin{proposition} \label{prop:semimartingale_ZN}
	 	Define $ \eta_N = \delta_N^{2-d} $ in the fixed radius case and $ \eta_N = 1 $ in the stable case.
	 	Then there exists a sequence of worthy martingale measures $ (M^N, N \geq 1) $ on $ \R_+ \times \R^d \times [0,1] $ such that, for any $ \phi \in E_1 $,
	 	\begin{align} \label{semimartingale_ZN}
	 	\langle Z^N_t, \phi \rangle = \int_{0}^{t} \langle Z^N_s, u \mathcal{L}^{N, \alpha} \phi - \mu \phi \rangle ds + M^N_t(\phi)
	 	\end{align}
	 	and the predictable variation process of $ M^N_t(\phi) $ is given by
	 	\begin{align} \label{qvar_MN}
	 	\left\langle M^N(\phi) \right\rangle_t = u^2 \int_{0}^{t} \langle \eta_N \Gamma^{\nu_\alpha^N}(\bm{\rho}^N_s), \phi \otimes \phi \rangle ds.
	 	\end{align}
	 \end{proposition}
 
     Note that the predictable part of $ \langle Z^N_t, \phi \rangle $ in \eqref{semimartingale_ZN} is linear in $ Z^N $.
     We can use this and the semigroup defined in \eqref{def_PNalpha} to write, for $ \phi \in E_1 $ \citep[Theorem~5.1]{walsh_introduction_1986},
     \begin{align} \label{stochastic_integral}
     \langle Z^N_t, \phi \rangle = \int_{[0,t] \times \R^d \times [0,1]} P^{N,\alpha}_{t-s}\phi(x,k) M^N(ds\, dx\, dk),
     \end{align}
     where the right hand side is defined as a stochastic integral against the martingale measure $ M^N $ (see Chapter~2 of \cite{walsh_introduction_1986}).
     This reduces the convergence of $ \langle Z^N_t, \phi \rangle $ to the convergence of a sequence of stochastic integrals.
     The main ingredients for this are the convergence of the sequence of martingale measures $ M^N $ and the convergence of the semigroups $ (P^{N,\alpha}_t, t \geq 0) $ as $ N \to \infty $.
     
     After proving Proposition~\ref{prop:semimartingale_ZN} in Subsection~\ref{subsec:semimartingale}, we shall prove that the sequence of semigroups $ (P^{N,\alpha}_t, t \geq 0) $ converges as $ N \to \infty $ to a semigroup $ (P^{(\alpha)}_t, t \geq 0) $ satisfying
     \begin{align*}
     P^{(\alpha)}_t \phi (x,k) := e^{-\mu t} \, \E[x]{\phi(X^{\alpha}_t,k)},
     \end{align*}
     for any bounded and measurable $ \phi : \R^d \times [0,1] \to \R $, where $ (X^\alpha_t, t \geq 0) $ is either Brownian motion or fractional Brownian motion, according as $ \alpha = 2 $ or not.
     This convergence is proved for an appropriate norm in Subsection~\ref{subsec:test_functions} (Lemma~\ref{lemma:test_functions}).
     
     Using a result adapted from \cite{walsh_introduction_1986} (Theorem~\ref{thm:walsh} below), this will allow us to prove, in Subsection~\ref{subsec:det_limit}, that $ \lbrace (\langle Z^N_t, \phi \rangle, t \geq 0), N \geq 1 \rbrace $ is tight in $ D(\R_+, \R) $ for any $ \phi \in \mathcal{S}(\R^d \times [0,1]) $, as well as the convergence of $ (\bm{\rho}^N_t, t \geq 0) $ to $ \lambda $, in probability locally uniformly in time, as $ N \to \infty $.
     
     Finally, we shall see in Subsection~\ref{subsec:convergence_Mn} that the right-hand-side of \eqref{qvar_MN} converges to $ t\, \langle \mathcal{Q}_\alpha, \phi \otimes \phi \rangle $ as $ N \to \infty $, where $ \mathcal{Q}_2 := \mathcal{Q} $ is defined in \eqref{def_Q} and $ \mathcal{Q}_\alpha $ is defined in \eqref{def_Qalpha} for $ \alpha \neq 2 $.
     We then deduce from this (in Lemma~\ref{lemma:cvg_MN}) the convergence in distribution of the sequence of martingale measures $ (M^N, N \geq 1) $ to a continuous martingale measure $ M $ such that, for any $ \phi \in \mathcal{S}(\R^d \times [0,1]) $,
     \begin{align*}
     \left\langle M(\phi) \right\rangle_t = t\, \langle \mathcal{Q}_\alpha, \phi \otimes \phi \rangle.
     \end{align*}
     
     Using this, \eqref{stochastic_integral} and the convergence of the semigroups $ P^{N,\alpha} $, we will apply a result on the convergence of stochastic integrals (recalled in Appendix~\ref{sec:proof_walsh}) to prove the convergence in distribution of $ \langle Z^N_t, \phi \rangle $ to $ \langle Z_t, \phi \rangle $, where
     \begin{align*}
     \langle Z_t, \phi \rangle := \int_{[0,t] \times \R^d \times [0,1]} P^{(\alpha)}_{t-s} \phi(x,k) M(ds\, dx\, dk),
     \end{align*}
     for any $ t \geq 0 $ and $ \phi \in \mathcal{S}(\R^d \times [0,1]) $.
     This convergence will easily be generalised to vectors of the form
     \begin{align*}
     	\left( \langle Z^N_{t_1}, \phi_1 \rangle, \ldots, \langle Z^N_{t_p}, \phi_p \rangle \right),
     \end{align*}
     with $ t_1, \ldots, t_p \in \R_+ $ and $ \phi_1, \ldots, \phi_p \in \mathcal{S}(\R^d \times [0,1]) $.
     We shall then conclude using the following theorem, which can be found in \citep[Theorem~6.15]{walsh_introduction_1986}.
     
     \begin{theorem} \label{thm:cvg_S'}
     	Let $ (X^N, N \geq 1) $ be a sequence of processes with sample paths in $ D([0,T]$, $\mathcal{S}'(\R^d\times [0,1])) $.
     	Suppose that
     	\begin{enumerate}[i)]
     		\item for each $ \phi \in \mathcal{S}(\R^d \times [0,1]) $, $ \lbrace (\langle X^N_t, \phi \rangle, t \in [0,T]), N \geq 1 \rbrace $ is tight,
     		\item for each $ \phi_1, \ldots, \phi_p $ in $ \mathcal{S}(\R^d \times [0,1]) $ and $ t_1, \ldots, t_p \in [0,T] $, the distribution of
     		\begin{align*}
     		(\langle X^N_{t_1}, \phi_1 \rangle, \ldots, \langle X^N_{t_p}, \phi_p \rangle)
     		\end{align*}
     		converges weakly on $ \R^p $.
     	\end{enumerate}
     	Then there exists a process $ (X_t, t \in [0,T]) $, with sample paths in $ D([0,T], \mathcal{S}'(\R^d \times [0,1])) $ such that $ X^N $ converges in distribution to $ X $.
     \end{theorem}
 
 	This last step is detailed in Subsection~\ref{subsec:proof_clt}.
	 
	 \subsection{The SLFV with mutations as a semimartingale} \label{subsec:semimartingale}
	 
	 Before proving Proposition~\ref{prop:semimartingale_ZN}, let us first state the semimartingale decomposition of the SLFV with mutations of Definition~\ref{def:slfv}, before rescaling space and time.
	 Recall that $ \lambda $ was defined in \eqref{def:lambda} as the Lebesgue measure on $ \R^d \times [0,1] $.
	 
	 \begin{proposition} \label{prop:slfv_semimartingale}
	 	Let $ (\rho_t, t \geq 0) $ be the SLFV with mutations of Definition~\ref{def:slfv} started from some $ \rho_0 \in \Xi $.
	 	Let $ (\F_t, t \geq 0) $ denote its natural filtration.
	 	For any $ \phi \in E_1 $,
	 	\begin{align} \label{martingale_slfv}
	 		\langle \rho_t, \phi \rangle - \langle \rho_0, \phi \rangle - \int_{0}^{t} \left\lbrace \mu \langle \lambda - \rho_s, \phi \rangle + u \int_{0}^{\infty} V_r \left\langle \rho_s, \overline{\overline{\phi}}(\cdot, r) - \phi \right\rangle \nu(dr) \right\rbrace ds
	 	\end{align}
	 	defines a (mean-zero) square integrable $ \F_t $-martingale with predictable variation process
	 	\begin{align*}
	 		u^2 \int_{0}^{t} \langle \Gamma^\nu(\rho_s), \phi \otimes \phi \rangle ds.
	 	\end{align*}
	 \end{proposition}
 	 
 	 The different terms appearing in this decomposition each correspond to a distinct evolutionary force.
 	 The term $ \mu \langle \lambda - \rho_s, \phi \rangle $ results from the mutations, the second term inside the integral in \eqref{martingale_slfv} is the spatial mixing resulting from the reproduction events (\textit{i.e.} the migration term) and the martingale part captures the fluctuations due to genetic drift, that is to say the randomness due to reproduction in a (locally) finite population.
 	 
 	 \begin{proof}[Proof of Proposition~\ref{prop:slfv_semimartingale}]
 	 	From Definition~\ref{def:slfv}, we have
 	 	\begin{multline*}
 	 		\lim_{\delta t \downarrow 0} \frac{1}{\delta t} \E{ \langle \rho_{t + \delta t}, \phi \rangle - \langle \rho_t, \phi \rangle }{ \rho_t = \rho } \\
 	 		\begin{aligned}
 	 		& = \mu \langle \lambda - \rho, \phi \rangle + \int_{\R^d} \int_{0}^{\infty} \frac{1}{V_r} \int_{B(x,r) \times [0,1]} \left\langle u \1_{B(x,r)} (\delta_{k_0} - \rho), \phi \right\rangle \rho(y, dk_0) dy \nu(dr) dx \\
 	 		& = \mu \langle \lambda - \rho, \phi \rangle + u \int_{0}^{\infty} V_r \left\langle \rho, \overline{\overline{\phi}}(\cdot,r) - \phi \right\rangle \nu(dr).
 	 		\end{aligned}
 	 	\end{multline*}
 	 	It follows (see \textit{e.g.} \citep[Proposition~4.1.7]{ethier_markov_1986}) that \eqref{martingale_slfv} defines a martingale.
 	 	To compute its variation process, write
 	 	\begin{multline*}
 	 		\lim_{\delta t \downarrow 0} \frac{1}{\delta t} \E{ \left( \langle \rho_{t + \delta t}, \phi\rangle - \langle \rho_t, \phi \rangle \right)^2 }{ \rho_t = \rho } \\ = \int_{\R^d} \int_{0}^{\infty} \frac{1}{V_r} \int_{B(x,r) \times [0,1]} \int_{(\R^d \times [0,1])^2} \phi(x_1, k_1) \phi(x_2, k_2) u^2 \1_{\lbrace \|x_1 - x\| < r \rbrace} \1_{\lbrace \|x_2 - x\| < r \rbrace} \\ \times ( \delta_{k_0}(dk_1) - \rho(x_1, dk_1) ) ( \delta_{k_0}(dk_2) - \rho(x_2, dk_2) ) dx_1 dx_2 \rho(y,dk_0) dy \nu(dr) dx.
 	 	\end{multline*}
 	 	Rearranging the integrals with respect to $ k_0 $, $ k_1 $ and $ k_2 $, this becomes
 	 	\begin{multline*}
 	 		\lim_{\delta t \downarrow 0} \frac{1}{\delta t} \E{ \left( \langle \rho_{t + \delta t}, \phi \rangle - \langle \rho_t, \phi \rangle \right)^2 }{ \rho_t = \rho } \\ = \int_{\R^d} \int_{0}^{\infty} \frac{1}{V_r} \int_{B(x,r)} \int_{(\R^d)^2} u^2 \1_{\lbrace \|x_1 - x\| < r \rbrace} \1_{\lbrace \|x_2 - x\| < r \rbrace} \left[ \int_{[0,1]} \phi(x_1, k_0) \phi(x_2, k_0) \rho(y, dk_0) \right. \\ - \int_{[0,1]} \phi(x_1, k_0) \rho(y, dk_0) \int_{[0,1]} \phi(x_2, k_2) \rho(x_2, dk_2) \\ - \int_{[0,1]} \phi(x_1, k_1) \rho(x_1, dk_1) \int_{[0,1]} \phi(x_2, k_0) \rho(y, dk_0) \\ + \left. \int_{[0,1]} \phi(x_1, k_1) \rho(x_1, dk_1) \int_{[0,1]} \phi(x_2, k_2) \rho(x_2, dk_2) \right] dx_1 dx_2 dy \nu(dr) dx.
 	 	\end{multline*}
 	 	By the definition of $ [\rho]_r $ and $ \Gamma^\nu(\rho) $ in \eqref{def_bracket_rho} and \eqref{def_Gamma}, this is
 	 	\begin{align*}
 	 		\lim_{\delta t \downarrow 0} \frac{1}{\delta t} \E{ \left( \langle \rho_{t + \delta t}, \phi \rangle - \langle \rho_t, \phi \rangle \right)^2 }{ \rho_t = \rho } = u^2 \langle \Gamma^\nu(\rho), \phi \otimes \phi \rangle.
 	 	\end{align*}
 	 	This concludes the proof of Proposition~\ref{prop:slfv_semimartingale}.
 	 \end{proof}
	 
	 We now use Proposition~\ref{prop:slfv_semimartingale} to prove Proposition~\ref{prop:semimartingale_ZN}.

	 \begin{proof}[Proof of Proposition~\ref{prop:semimartingale_ZN}]
	     Recall that we have set
	 \begin{align*}
	 	\bm{\rho}^N_t(x, dk) = \rho^N_{Nt/\delta_N^\alpha}(x/\delta_N, dk)
	 \end{align*}
	 where $ \alpha = 2 $ in the fixed radius case.
	 For $ \phi \in E_1 $,
	 \begin{align} \label{rescaling_slfv}
	 	\langle \bm{\rho}^N_t, \phi \rangle = \langle \rho^N_{Nt/\delta_N^\alpha}, \phi_N \rangle
	 \end{align}
	 with
	 \begin{align} \label{def_phiN}
	 	\phi_N(x, k) = \delta_N^d \phi(\delta_N x, k).
	 \end{align}
	 Let $ (\mathcal{M}^N_t(\phi), t \geq 0) $ denote the martingale defined by \eqref{martingale_slfv} in Proposition~\ref{prop:slfv_semimartingale}, \textit{i.e.}
	 \begin{align*}
	 	\mathcal{M}^N_t(\phi) = \langle \rho^N_t, \phi \rangle - \langle \rho^N_0, \phi \rangle - \int_{0}^{t} \left\lbrace \mu_N \langle \lambda - \rho^N_s, \phi \rangle + u_N \int_{0}^{\infty} V_r \left\langle \rho^N_s, \overline{\overline{\phi}}(\cdot, r) - \phi \right\rangle \nu_\alpha(dr) \right\rbrace ds,
	 \end{align*}
	 where we recall that $ \nu_\alpha = \delta_R $ for $ \alpha = 2 $.
	 Then, by \eqref{rescaling_slfv},
	 \begin{multline*}
	 	\langle \bm{\rho}^N_t, \phi \rangle = \langle \rho^N_0, \phi_N \rangle + \int_{0}^{Nt / \delta_N^\alpha} \left\lbrace \mu_N \langle \lambda - \rho^N_s, \phi_N \rangle + u_N \int_{0}^{\infty} V_r \left\langle \rho^N_s, \overline{\overline{\phi_N}}(\cdot, r) - \phi_N \right\rangle \nu_\alpha(dr) \right\rbrace ds \\ + \mathcal{M}^N_{Nt / \delta_N^\alpha}(\phi_N).
	 \end{multline*}
	 But, by a simple change of variables,
	 \begin{align*}
	 	\overline{\phi_N}(x,k,r) = \delta_N^d \overline{\phi}(\delta_N x, k, \delta_N r).
	 \end{align*}
	 As a result, replacing $ \mu_N = \delta_N^\alpha \frac{\mu}{N} $ and $ u_N = \frac{u}{N} $ and changing variables in the time integral,
	 \begin{align} \label{semimartingale_rescaled}
	 	\langle \bm{\rho}^N_t, \phi \rangle = \langle \bm{\rho}^N_0, \phi \rangle + \int_{0}^{t} \left\lbrace \mu \langle \lambda - \bm{\rho}^N_s, \phi \rangle + u \langle \bm{\rho}^N_s, \mathcal{L}^{N,\alpha} \phi \rangle \right\rbrace ds + \mathcal{M}^N_{Nt / \delta_N^\alpha} (\phi_N)
	 \end{align}
	 recalling the definition of $ \mathcal{L}^{N,\alpha} $ in \eqref{def_LN}.
	 Recall also that $ Z^N_t $ was defined by
	 \begin{align*}
	 	Z^N_t = (N \eta_N)^{1/2} (\bm{\rho}^N_t - \lambda).
	 \end{align*}
	 Subtracting $ \langle \lambda, \phi \rangle $ on both sides of \eqref{semimartingale_rescaled} and multiplying by $ (N \eta_N)^{1/2} $, we see that $ Z^N $ satisfies \eqref{semimartingale_ZN} where the martingale measure $ M^N $ is defined by
	 \begin{align*}
	 M^N_t(\phi) = (N \eta_N)^{1/2} \mathcal{M}^N_{Nt / \delta_N^\alpha}(\phi_N).
	 \end{align*}
	 Then, by Proposition~\ref{prop:slfv_semimartingale},
	 \begin{align*}
	 	\left\langle M^N(\phi) \right\rangle_t = (N \eta_N)\, u_N^2 \int_{0}^{Nt / \delta_N^\alpha} \left\langle \Gamma^{\nu_\alpha}(\rho^N_s), \phi_N \otimes \phi_N \right\rangle ds.
	 \end{align*}
	 Again, by a change of variables,
	 \begin{align*}
	 	\langle \Gamma^{\nu_\alpha}(\rho^N_{Ns/\delta_N^\alpha}), \phi_N \otimes \phi_N \rangle = \delta_N^\alpha \langle \Gamma^{\nu_\alpha^N}(\bm{\rho}^N_s), \phi \otimes \phi \rangle
	 \end{align*}
	 where $ \nu_\alpha^N $ was defined in \eqref{def_nu_N_alpha}.
	 As a result,
	 \begin{align} \label{qvar_rescaled}
	 	\left\langle M^N(\phi) \right\rangle_t = u^2 \int_{0}^{t} \langle \eta_N \Gamma^{\nu_\alpha^N}(\bm{\rho}^N_s), \phi \otimes \phi \rangle ds.
	 \end{align}
	 Together, \eqref{semimartingale_rescaled} and \eqref{qvar_rescaled} yield \eqref{semimartingale_ZN} and \eqref{qvar_MN}.
	 
 	 	We are left with proving that $ M^N $ is worthy (see the definition in Chapter~2 of \cite{walsh_introduction_1986}).
 	 	To do this, define $ \abs{\Gamma}^{\nu}(\rho) $ by
 	 	\begin{multline*}
 	 		\abs{\Gamma}^{\nu}(\rho)(x_1, x_2, dk_1 dk_2) = \int_{0}^{\infty} V_{2,r}(x_1, x_2) \Big[ [\rho]_r(x_1, x_2, dk_1) \delta_{k_1}(dk_2) + \rho(x_1, dk_1) [\rho]_r(x_1, x_2, dk_2) \\ + [\rho]_r(x_1, x_2, dk_1) \rho(x_2, dk_2) + \rho(x_1, dk_1) \rho(x_2, dk_2) \Big] \nu(dr).
 	 	\end{multline*}
 	 	Then the measure
 	 	\begin{align} \label{def_KN}
 	 		K_N(dt dx_1 dk_1 dx_2 dk_2) = u^2 \eta_N \abs{\Gamma}^{\nu_\alpha^N}(\bm{\rho}^N_{t^-})(x_1, x_2, dk_1 dk_2) dx_1 dx_2 dt
 	 	\end{align}
 	 	is positive definite and symmetric in $ (x_1, k_1) $, $ (x_2, k_2) $.
 	 	In addition, for any $ A, B \subset \R^d \times [0,1] $, $ (K_N([0,t] \times A \times B), t \geq 0) $ is predictable and for any rectangle $ \Lambda \subset [0,\infty) \times (\R^d \times [0,1])^2 $,
 	 	\begin{align*}
 	 		\abs{Q_N(\Lambda)} \leq K_N(\Lambda), \quad \text{ a.s.}
 	 	\end{align*}
 	 	where $ Q_N $ is the covariation measure of $ M^N $, see \eqref{qvar_rescaled}.
 	 	Thus, $ K_N $ is a dominating measure for $ M^N $, and $ M^N $ is a worthy martingale measure.
 	 	This concludes the proof of Proposition~\ref{prop:semimartingale_ZN}.
 	 \end{proof}

  \subsection{Bound on the dominating measures} \label{subsec:bound_qvar}
  
  We now prove the following.
  
  \begin{lemma} \label{lemma:bound_qvar_MN}
  	For $ N \geq 1 $, let $ K_N $ be the measure defined in \eqref{def_KN}.
  	There exists a constant $ \newCst{KN} > 0 $ such that, for all $ N \geq 1 $, for all $ 0 \leq s \leq t $ and for all $ \phi \in \mathcal{S}(\R^d \times [0,1]) $,
  	\begin{align*}
  	\langle K_N, \1_{[s,t]} \phi \otimes \psi \rangle \leq \Cst{KN} \abs{t-s} \left( \| \phi \|_1 \left\| \psi \right\|_1 + \| \phi \|_2 \left\| \psi \right\|_2 \right).
  	\end{align*}
  \end{lemma}
  
  \begin{proof}[Proof of Lemma~\ref{lemma:bound_qvar_MN}]
  	From \eqref{def_KN},
  	\begin{align*}
  	\langle K_N, \1_{[s,t]} \phi \otimes \psi \rangle = u^2 \int_{s}^{t} \langle \eta_N \abs{\Gamma}^{\nu_\alpha^N}(\bm{\rho}^N_v), \phi \otimes \psi \rangle dv.
  	\end{align*}
  	This is bounded by
  	\begin{align} \label{bound_qvar_1}
  	4 u^2 \abs{t-s} \int_{0}^{\infty} \int_{(\R^d)^2} V_{\delta_N r}(x_1, x_2) \sup_{k \in [0,1]} \abs{\phi(x_1,k)} \sup_{k \in [0,1]} \abs{\psi(x_2,k)} dx_1 dx_2 \eta_N \delta_N^{-(d+\alpha)} \nu_\alpha(dr).
  	\end{align}
  	We split the integral over $ r $ on $ [0,1/\delta_N] $ and $ (1/\delta_N, \infty) $.
  	In the first integral, we use the Cauchy-Schwarz inequality and the fact that
  	\begin{align*}
  	\int_{\R^d} V_{2,r}(x_1, x_2) dx_2 = V_r^2,
  	\end{align*}
  	to obtain
  	\begin{align*}
  	\int_{(\R^d)^2} V_{\delta_N r}(x_1, x_2) \sup_{k \in [0,1]} \abs{\phi(x_1,k)} \sup_{k \in [0,1]} \abs{\phi(x_2,k)} dx_1 dx_2 \leq V_{\delta_N r}^2 \left\| \phi \right\|_2 \left\| \psi \right\|_2.
  	\end{align*}
  	In the second integral, we simply use $ V_{2,r}(x_1, x_2) \leq V_r $.
  	As a result, \eqref{bound_qvar_1} is bounded by
  	\begin{align*}
  	4 u^2 \abs{t-s} \left\lbrace \int_{0}^{1/\delta_N} V_{\delta_N r}^2 \eta_N \delta_N^{-(d+\alpha)} \nu_\alpha(dr) \| \phi \|_2 \left\| \psi \right\|_2 + \int_{1/\delta_N}^{\infty} V_{\delta_N r} \eta_N \delta_N^{-(d+\alpha)} \nu_\alpha(dr) \| \phi \|_1 \left\| \psi \right\|_1 \right\rbrace.
  	\end{align*}
  	In the fixed radius case, clearly the second integral vanishes for $ N $ large enough and the first one is
  	\begin{align*}
  	V_{\delta_N R}^2 \delta_N^{2-d} \delta_N^{-(d+2)} = V_R^2 < \infty.
  	\end{align*}
  	In the stable case, $ \eta_N = 1 $ and
  	\begin{align*}
  	\int_{0}^{1/\delta_N} V_{\delta_N r}^2 \delta_N^{-(d+\alpha)} \nu_\alpha(dr) &= \int_{1}^{1/\delta_N} V_{\delta_N r}^2 \frac{dr}{r(\delta_N r)^{\alpha+d}} \\
  	&= \int_{\delta_N}^{1} V_r^2 \frac{dr}{r^{1+\alpha+d}} \\
  	&\leq \int_{0}^{1} V_r^2 \frac{dr}{r^{1+\alpha+d}} < \infty
  	\end{align*}
  	since $ d > \alpha $.
  	For the second integral,
  	\begin{align*}
  	\int_{1/\delta_N}^{\infty} V_{\delta_N r} \delta_N^{-(d+\alpha)} \nu_\alpha(dr) = \int_{1}^{\infty} V_r \frac{dr}{r^{1+\alpha+d}} < \infty.
  	\end{align*}
  	The statement of Lemma~\ref{lemma:bound_qvar_MN} then follows.
  \end{proof}

\subsection[Convergence of the semigroups P]{Convergence of the semigroups $ P^{N,\alpha} $} \label{subsec:test_functions}

Let us begin this subsection by the following observations.
Recall the definition of $ \mathcal{L}^{N,\alpha} $ in \eqref{def_LN}.
In the fixed radius case, substituting $ \nu_\alpha(dr) = \delta_R(dr) $ in \eqref{def_LN},
\begin{align*}
\mathcal{L}^{N, 2} \phi (x, k) = \frac{V_R}{\delta_N^2} \left( \overline{\overline{\phi}}(x, k, \delta_N R) - \phi(x,k) \right).
\end{align*}
Recalling the notation \eqref{def_average} and writing a Taylor expansion inside the spatial average, we obtain (see Proposition~\ref{prop:averages})
\begin{align} \label{cvg_L_alpha_laplace}
\| \mathcal{L}^{N,\alpha} \phi - \mathcal{D}_2 \phi \|_{q} \leq V_R R^4 \frac{d^3}{3} (\delta_N)^2 \max_{\abs{\beta} = 4} \| \partial_\beta \phi \|_q,
\end{align}
for any $ \phi : \R^d \times [0,1] \to \R $ which admits continuous and $ \| \cdot \|_q $-bounded spatial derivatives of order up to four, where $ \mathcal{D}_2 $ is defined as
\begin{align} \label{def_D2}
\mathcal{D}_2 \phi := V_R \frac{R^2}{d+2} \Delta \phi,
\end{align}
On the other hand, in the stable case,
\begin{align*}
\mathcal{L}^{N, \alpha} \phi(x,k) &= \delta_N^{-\alpha} \int_{1}^{\infty} V_r \left( \overline{\overline{\phi}}(x,k,\delta_N r) - \phi(x,k) \right) \frac{dr}{r^{1+\alpha+d}} \\
&= \int_{\delta_N}^{\infty} V_r \left( \overline{\overline{\phi}}(x,k,r) - \phi(x,k) \right) \frac{dr}{r^{1+\alpha+d}} \\
&= \int_{\R^d} \Phi_\alpha^{\delta_N}(\|x-y\|) (\phi(y, k) - \phi(x,k)) dy,
\end{align*}
with
\begin{align*}
\Phi_\alpha^\delta(\|x-y\|) = \int_{\delta}^{\infty} \frac{V_{2,r}(x,y)}{V_r} \frac{dr}{r^{1+\alpha+d}}.
\end{align*}
Hence, from Proposition~\ref{prop:fractional_laplacian},
\begin{align} \label{cvg_L_alpha_frac}
\| \mathcal{L}^{N,\alpha} \phi - \mathcal{D}_\alpha \phi \|_q \leq \Cst{fracL} (\delta_N)^{2-\alpha} \max_{\abs{\beta} = 2} \| \partial_\beta \phi \|_q,
\end{align}
for any $ \phi : \R^d \times [0,1] \to \R $ twice continuously differentiable in the space variable and such that the right hand side is finite.

Moreover, the operator $ u \mathcal{D}_\alpha $ (or $ u \mathcal{D}_2 $ in the fixed radius case), defined on the set of twice continuously differentiable functions, is the generator of a Markov process on $ \R^d $, denoted by $ (X^{\alpha}_t, t \geq 0) $, which is Brownian motion in the fixed radius case and fractional Brownian motion in the stable case.
We thus let $ (P^{(\alpha)}_t, t \geq 0) $ denote the strongly continuous semigroup acting on the space of bounded and measurable functions on $ \R^d \times [0,1] $ generated by $ u \mathcal{D}_\alpha - \mu $, \textit{i.e.}
\begin{align*}
P^{(\alpha)}_t \phi (x,k) := e^{-\mu t} \, \E[x]{\phi(X^{\alpha}_t,k)},
\end{align*}
for any bounded and measurable $ \phi : \R^d \times [0,1] \to \R $ and where $ \E[x]{\cdot} $ denotes the expectation with respect to the distribution of $ (X^{\alpha}_t, t \geq 0) $ started from $ X^{\alpha}_0 = x $.

Using this, we now prove the following lemma.

\begin{lemma}[Convergence of the semigroups $ P^{N,\alpha} $] \label{lemma:test_functions}
	For any $ q \geq 1 $, $ \phi \in E_q $ and for all $ t \geq 0 $, $ N \geq 1 $,
	\begin{align} \label{bound_test_functions}
	\| P^{N,\alpha}_t \phi \|_q \leq e^{-\mu t} \| \phi \|_q.
	\end{align}
	Furthermore for any multi-index $ \beta \in \N^d $ and $ \phi \in \mathcal{S}(\R^d \times [0,1]) $,
	\begin{align} \label{bound_deriv_test_functions}
	\| \partial_\beta P^{N,\alpha}_t \phi \|_q \leq e^{-\mu t} \| \partial_\beta \phi \|_q.
	\end{align}
	Finally, there exists a constant $ \newCst{cvgpsi} > 0 $ such that for all $ N \geq 1 $ and $ t \geq 0 $,
	\begin{align} \label{convergence_test_functions}
	\| P^{N,\alpha}_t \phi - P^{(\alpha)}_t \phi \|_q \leq \Cst{cvgpsi} t e^{-\mu t} (\delta_N)^\gamma  \max_{0 \leq \abs{\beta} \leq 4} \| \partial_\beta \phi \|_q
	\end{align}
	where $ \gamma = 2 $ in the fixed radius case and $ \gamma = 2 -\alpha $ in the stable case.
\end{lemma}

\begin{proof}[Proof of Lemma~\ref{lemma:test_functions}]
	Recall \eqref{def_PNalpha} and the definition of $ (X^{N,\alpha}_t, t \geq 0) $.
	Since the jumps of $ X^{N,\alpha} $ do not depend on its initial position,
	\begin{align} \label{translation_invariance_1}
	\E[x]{\phi(X^{N,\alpha}_t, k)} = \E[0]{ \phi(x + X^{N,\alpha}_t, k) },
	\end{align}
	and so, by Fubini's theorem,
	\begin{align} \label{translation_invariance}
	\int_{\R^d} \E[x]{ \phi(X^{N,\alpha}_t, k) } dx = \int_{\R^d} \phi(x, k) dx.
	\end{align}
	Moreover, for $ q \geq 1 $, using Jensen's inequality in the second line,
	\begin{align*}
	\| P^{N,\alpha}_t \phi \|_q &= \left( \int_{\R^d} \sup_{k \in [0,1]} \abs{ e^{-\mu t} \E[x]{ \phi(X^{N,\alpha}, k) } }^q dx \right)^{1/q} \\
	&\leq e^{-\mu t} \left( \int_{\R^d} \E[x]{\sup_{k \in [0,1]} \abs{\phi(X^{N,\alpha}_t, k)}^q } dx \right)^{1/q}.
	\end{align*}
	Then, applying \eqref{translation_invariance} with $ \phi(\cdot,k) $ replaced by $ \sup_{k \in [0,1]} \abs{\phi(\cdot,k)}^q $, we obtain
	\begin{align*}
		\| P^{N,\alpha}_t \phi \|_q \leq e^{-\mu t} \| \phi \|_q.
	\end{align*}
	This proves the first part of the statement of Lemma~\ref{lemma:test_functions}.
	Using \eqref{translation_invariance_1} again, we see that
	\begin{align*}
	\partial_\beta P^{N,\alpha}_t \phi(x,k) &= e^{-\mu t} \E[x]{ \partial_\beta \phi(X^{N,\alpha}_t, k) } \\
	&= P^{N,\alpha}_t \partial_\beta \phi (x,k).
	\end{align*}
	Thus \eqref{bound_deriv_test_functions} follows from \eqref{bound_test_functions} applied to $ \partial_\beta \phi $.
	To prove \eqref{convergence_test_functions}, we note that the above inequalities also apply to the semigroup $ P^{(\alpha)}_t $, replacing $ X^{N,\alpha}_t $ by $ X^\alpha_t $, and we let
	\begin{align*}
	\psi_t = P^{N,\alpha}_t \phi - P^{(\alpha)}_t \phi.
	\end{align*}
	Then
	\begin{align*}
	\partial_t \psi_t = (u\mathcal{L}^{N,\alpha}-\mu) \psi_t + u(\mathcal{L}^{N,\alpha} - \mathcal{D}_\alpha) P^{(\alpha)}_t \phi.
	\end{align*}
	Together with $ \psi_0 = 0 $, this implies
	\begin{align*}
	\psi_t = u \int_{0}^{t} P^{N,\alpha}_{t-s} \left( \mathcal{L}^{N,\alpha} - \mathcal{D}_\alpha \right) P^{(\alpha)}_s \phi \, ds.
	\end{align*}
	Using \eqref{bound_test_functions} and the triangle inequality, we obtain
	\begin{align*}
	\left\| \psi_t \right\|_q \leq u \int_{0}^{t} e^{-\mu(t-s)} \left\| (\mathcal{L}^{N,\alpha} - \mathcal{D}_\alpha) P^{(\alpha)}_s \phi \right\|_q ds.
	\end{align*}
	Using \eqref{cvg_L_alpha_laplace} in the fixed radius case and \eqref{cvg_L_alpha_frac} in the stable case, we see that there exists a constant $ \newCst{psit} > 0 $ such that
	\begin{align*}
	\left\| \psi_t \right\|_q \leq u \int_{0}^{t} e^{-\mu(t-s)} \Cst{psit} (\delta_N)^\gamma \max_{\abs{\beta} \leq 4} \left\| \partial_\beta P^{(\alpha)}_s \phi \right\|_q ds
	\end{align*}
	where $ \gamma = 2 $ in the fixed radius case and $ \gamma = 2-\alpha $ in the stable case.
	Finally, using \eqref{bound_deriv_test_functions} applied to $ P^{(\alpha)}_s $, we obtain
	\begin{align*}
	\left\| \psi_t \right\|_q &\leq u \Cst{psit} (\delta_N)^\gamma \int_{0}^{t} e^{-\mu(t-s)} e^{-\mu s} \max_{\abs{\beta} \leq 4} \left\| \partial_\beta \phi \right\|_q ds \\
	&\leq u \Cst{psit} t e^{-\mu t} (\delta_N)^\gamma \max_{\abs{\beta} \leq 4} \left\| \partial_\beta \phi \right\|_q.
	\end{align*}
	This concludes the proof of Lemma~\ref{lemma:test_functions}.
\end{proof}

    \subsection{Convergence of the rescaled SLFV to a deterministic limit} \label{subsec:det_limit}
    
  With the above results, it becomes possible to prove the first part of the two central limit theorems, \textit{i.e.} the convergence of $ (\bm{\rho}^N_t, t \in [0,T]) $ to $ \lambda $.
  
  \begin{proposition}[Convergence to the deterministic limit] \label{prop:convergence_deterministic_limit}
  	Under the assumptions of either Theorem~\ref{thm:clt_fixed_radius} or Theorem~\ref{thm:clt_stable}, for any $ T > 0 $,
  	\begin{align*}
  	\lim_{N \to \infty} \E{ \sup_{t \in [0,T]} d(\bm{\rho}^N_t, \lambda) } = 0.
  	\end{align*}
  \end{proposition}
  
  To prove this, we apply the following result, which is adapted from \citep[Theorem~7.13]{walsh_introduction_1986}.
  We detail its proof in Appendix~\ref{sec:proof_walsh}.
  
  \begin{theorem}[Adapted from Theorem~7.13 in \cite{walsh_introduction_1986}] \label{thm:walsh}
  	Let $ (M^N, N \geq 1) $ be a sequence of worthy martingale measures on $ \R^d \times [0,1] $ with dominating measures $ (K_N, N \geq 1) $ such that there exist $ \Cst{KN} > 0 $ and $ k \geq 1 $ with, for any $ 0 \leq s \leq t $, $ \phi \in \mathcal{S}(\R^d \times [0,1]) $ and for all $ N \geq 1 $,
  	\begin{align} \label{condition_K_N}
  	\int_{[s,t] \times (\R^d \times [0,1])^2} \phi(x_1, k_1) \phi(x_2, k_2) K_N(dr dx_1 dk_1 dx_2 dk_2) \leq \Cst{KN} \abs{t-s} \sum_{q=1}^{k} \left\| \phi \right\|_q^2,
  	\end{align}
  	almost surely.
  	Let $ \psi^N : \lbrace (s, t) : 0 \leq s \leq t \rbrace \times \R^d \times [0,1] \to \R $ be a sequence of deterministic functions such that:
  	\begin{enumerate}[i)]
  		\item \label{psi_schwartz} for any $ 0 \leq s \leq t $, $ N \geq 1 $, $ \psi^N_{s,t} \in \mathcal{S}(\R^d \times [0,1]) $,
  		\item \label{psi_continuous} both $ t \mapsto \psi^N_{s,t} $ and $ s \mapsto \psi^N_{s,t} $ are continuous,
  		\item \label{psi_bounds} there exist $ \newCst{psiN} > 0 $, $ \newCst{psiNtime} > 0 $ and $ \mu > 0 $ such that, for all $ N \geq 1 $, for all $ 0 \leq s \leq t $ and $ q \in [1, k] $,
  		\begin{align} \label{condition_psi_1}
  		\left\| \psi^N_{s,t} \right\|_q \leq \Cst{psiN} e^{-\mu(t-s)}
  		\end{align}
  		and, for all $ 0 \leq s' \leq s \leq t \leq t' $,
  		\begin{align} \label{condition_psi_2}
  		\left\| \psi^N_{s,t'} - \psi^N_{s,t} \right\|_q \leq \Cst{psiNtime} \abs{t'-t} e^{-\mu(t-s)}, && \left\| \psi^N_{s',t} - \psi^N_{s,t} \right\|_q \leq \Cst{psiNtime} \abs{s'-s} e^{-\mu(t-s)},
  		\end{align}
  	\end{enumerate}
  	Then the sequence of real-valued processes $ (U^N_t, t \geq 0) $ defined by
  	\begin{align*}
  	U^N_t = \int_{[0,t] \times \R^d \times [0,1]} \psi^N_{s,t}(x,k) M^N(ds\, dx\, dk)
  	\end{align*}
  	satisfies, for all $ T > 0 $, $ N \geq 1 $,
  	\begin{align} \label{sup_bound_U}
  	\E{ \sup_{t \in [0,T]} \abs{ U^N_t }^2 } \leq C_T\, \Cst{KN}\, k \left( \Cst{psiN}^2 + \Cst{psiNtime}^2 \right),
  	\end{align}
  	where the constant $ C_T $ only depends on $ T $.
  	In addition, the sequence $ (U^N, N \geq 1) $ is tight in $ D(\R_+, \R) $.
  \end{theorem}

  Using this, we can prove Proposition~\ref{prop:convergence_deterministic_limit}.
  
  \begin{proof}[Proof of Proposition~\ref{prop:convergence_deterministic_limit}]
  	We first check that the assumptions of Theorem~\ref{thm:walsh} are satisfied with
  	\begin{align*}
  	\psi^N_{s,t} = P^{N,\alpha}_{t-s} \phi.
  	\end{align*}
  	Clearly, \eqref{condition_K_N} (with $ k = 2 $) follows from Lemma~\ref{lemma:bound_qvar_MN}.
  	Conditions \ref{psi_schwartz} and \ref{psi_continuous} are also clearly satisfied.
  	In addition, \eqref{condition_psi_1} follows from \eqref{bound_test_functions} with
  	\begin{align} \label{Cst_psi_N}
  		\Cst{psiN} = \max_{q \in \lbrace 1, 2 \rbrace} \left\| \phi \right\|_q.
  	\end{align}
  	To prove \eqref{condition_psi_2}, write, for $ t \geq 0 $,
  	\begin{align*}
  	P^{N,\alpha}_t \phi - \phi = \int_{0}^{t} P^{N,\alpha}_s (u \mathcal{L}^{N,\alpha} - \mu) \phi\, ds.
  	\end{align*}
  	By the triangle inequality and \eqref{bound_test_functions},
  	\begin{align*}
  	\| P^{N,\alpha}_t \phi - \phi \|_q &\leq \int_{0}^{t} \| P^{N,\alpha}_s (u \mathcal{L}^{N,\alpha} - \mu) \phi \|_q ds \\
  	&\leq \int_{0}^{t} e^{-\mu s} \| (u \mathcal{L}^{N,\alpha} - \mu) \phi \|_q ds.
  	\end{align*}
  	Using Proposition~\ref{prop:averages} in the fixed radius case and Proposition~\ref{prop:fractional_laplacian} in the stable case,
  	\begin{align*}
  	\| (u \mathcal{L}^{N,\alpha} - \mu) \phi \|_q \leq C \max_{0 \leq \abs{\beta} \leq 2} \| \partial_\beta \phi \|_q,
  	\end{align*}
  	for some $ C > 0 $.
  	As a result, there exists a constant $ \newCst{continuitypsi} > 0 $ such that
  	\begin{align*}
  	\| P^{N,\alpha}_t \phi - \phi \|_q \leq \Cst{continuitypsi}\, t \max_{0\leq \abs{\beta} \leq 2} \| \partial_\beta \phi \|_q. 
  	\end{align*}
  	Then, for $ 0 \leq s \leq t \leq t' $,
  	\begin{align*}
  	\| \psi^N_{s,t'} - \psi_{s,t}^N \|_q &= \| P^{N,\alpha}_{t'-s}\phi - P^{N,\alpha}_{t-s}\phi \|_q \\
  	&\leq \Cst{continuitypsi} \abs{t'-t} \max_{0\leq \abs{\beta} \leq 2} \| \partial_\beta P^{N,\alpha}_{t-s} \phi \|_q \\
  	&\leq \Cst{continuitypsi} \abs{t'-t} e^{-\mu(t-s)} \max_{0\leq \abs{\beta} \leq 2} \| \partial_\beta \phi \|_q,
  	\end{align*}
  	where we have used \eqref{bound_deriv_test_functions} in the last line.
  	This proves the first part of \eqref{condition_psi_2} with
  	\begin{align} \label{Cst_psi_time}
  		\Cst{psiNtime} = \Cst{continuitypsi} \max_{q \in \lbrace 1, 2 \rbrace} \max_{0\leq \abs{\beta} \leq 2} \left\| \partial_\beta \phi \right\|_q.
  	\end{align}
  	The second part is proved in exactly the same way.
  	We can thus combine Theorem~\ref{thm:walsh} with \eqref{Cst_psi_N} and \eqref{Cst_psi_time} to obtain
  	\begin{align*}
  	\E{ \sup_{t \in [0,T]} \abs{ \langle Z^N_t, \phi \rangle }^2 }^{1/2} &\leq \left( 2 C_T\, \Cst{KN}\, \left( \Cst{psiN}^2 + \Cst{psiNtime}^2 \right) \right)^{1/2} \\
  	&\leq \newCst{supZ} \max_{q \in \lbrace 1, 2 \rbrace} \max_{0 \leq \abs{\beta} \leq 2} \left\| \partial_\beta \phi \right\|_q,
  	\end{align*}
  	for some $ \Cst{supZ} > 0 $ for all $ N \geq 1 $, $ T > 0 $ and $ \phi \in \mathcal{S}(\R^d \times [0,1]) $.
  	Then by the definition of the metric $ d $ in \eqref{def:d} and that of $ Z^N $, we have
  	\begin{align*}
  	\E{ \sup_{t \in [0,T]} d(\bm{\rho}^N_t, \lambda) } \leq \frac{\Cst{supZ}}{(N \eta_N)^{1/2}} \sum_{n = 1}^{\infty} \frac{1}{2^n} \max_{q \in \lbrace 1, 2 \rbrace} \max_{0 \leq \abs{\beta} \leq 2} \left\| \partial_\beta \phi_n \right\|_q.
  	\end{align*}
  	But recall from \eqref{phi_n_bound} that $ \| \partial_\beta \phi_n \|_q \leq \Cst{phi_n} $ for all $ n \geq 1 $, $ q \in \lbrace 1, 2 \rbrace $ and $ \beta \in \N^d $ with $ 0 \leq \abs{\beta} \leq 2 $.
  	Hence
  	\begin{align*}
  	\E{ \sup_{t \in [0,T]} d(\bm{\rho}^N_t, \lambda) } \leq \frac{\Cst{supZ} \Cst{phi_n}}{(N \eta_N)^{1/2}},
  	\end{align*}
  	and Proposition~\ref{prop:convergence_deterministic_limit} is proved.
  \end{proof}

	Note that, in passing, we have also proved the following.
	
	\begin{lemma}[Tightness] \label{lemma:tightness}
		For any $ \phi \in \mathcal{S}(\R^d \times [0,1]) $, the sequence of processes $ (\langle Z^N_t, \phi \rangle)_{t \geq 0} $, $ N \geq 1 $ is tight in $ D(\R_+, \R) $.
	\end{lemma}

\subsection{Convergence of the martingale measures} \label{subsec:convergence_Mn}

The aim of this subsection is to prove the following result.

\begin{lemma}[Convergence of the martingale measures] \label{lemma:cvg_MN}
	The sequence of martingale measures $ (M^N, N \geq 1) $ converges in distribution in $ D(\R_+,\mathcal{S}'(\R^d \times [0,1])) $ to a continuous martingale measure $ M $ such that, for any $ \phi \in \mathcal{S}(\R^d \times [0,1]) $,
	\begin{align*}
	\left\langle M(\phi) \right\rangle_t = t\, \langle \mathcal{Q}_\alpha, \phi \otimes \phi \rangle.
	\end{align*}
\end{lemma}

Lemma~\ref{lemma:cvg_MN} will result from the following, which we prove below.

\begin{lemma} \label{lemma:convergence_qvar_MN}
	For any $ \phi \in \mathcal{S}(\R^d \times [0,1]) $,
	\begin{align} \label{bound_jumps_MN}
	\sup_{t \geq 0} \abs{M^N_t(\phi) - M^N_{t^-}(\phi)} \cvgas{N} 0 \qquad \text{ almost surely.}
	\end{align}
	Furthermore, for any $ t \geq 0 $,
	\begin{align} \label{convergence_qvar_MN}
	\left\langle M^N(\phi) \right\rangle_t \cvgas{N} t \left\langle Q_\alpha, \phi \otimes \phi \right\rangle \qquad \text{ in probability,}
	\end{align}
	where $ Q_\alpha $ is defined in \eqref{def_Qalpha} in the stable case and $ \mathcal{Q}_2 = \mathcal{Q} $ is defined in \eqref{def_Q} in the fixed radius case.
\end{lemma}

Let us show how this implies Lemma~\ref{lemma:cvg_MN}.

\begin{proof}[Proof of Lemma~\ref{lemma:cvg_MN}]
	Lemma~\ref{lemma:convergence_qvar_MN}, together with Theorem~\ref{thm:Jacod_Sh} (in Appendix~\ref{sec:Jacod-Sh}) implies that, for all $ \phi \in \mathcal{S}(\R^d \times [0,1]) $, $ (M^N_t(\phi), t \geq 0) $ converges to $ (M_t(\phi), t \geq 0) $ in distribution in $ D(\R_+, \R) $.
	The sequence $ (M^N, N \geq 1) $ is thus tight in $ D(\R_+, \mathcal{S}'(\R^d \times [0,1])) $ by Mitoma's theorem (Theorem~\ref{thm:mitoma} in Appendix~\ref{sec:walsh}).
	
	Also, by polarisation, we can recover $ \langle M^N(\phi_i), M^N(\phi_j) \rangle_t $ from $ \langle M^N(\phi_i+\phi_j) \rangle_t $ and $ \langle M^N(\phi_i-\phi_j) \rangle_t $, and Theorem~\ref{thm:Jacod_Sh} is also satisfied by vectors of the form
	\begin{align*}
		\left( M^N_t(\phi_1), \ldots, M^N_t(\phi_k) \right)_{t \geq 0}.
	\end{align*}
	As a result, $ M^N $ satisfies the assumptions of Theorem~\ref{thm:cvg_S'}, and the sequence $ (M^N, N \geq 1) $ converges in distribution to $ M $ in $ D(\R_+, \mathcal{S}'(\R^d \times [0,1])) $.
\end{proof}

We now turn to the proof of Lemma~\ref{lemma:convergence_qvar_MN}.

\begin{proof}[Proof of Lemma~\ref{lemma:convergence_qvar_MN}]
	We first find a bound on the jumps of $ (M^N_t(\phi), t \geq 0) $, \textit{i.e.} \eqref{bound_jumps_MN}.
	By the definition of $ M^N $ in Proposition~\ref{prop:semimartingale_ZN},
	\begin{align*}
	M^N_t(\phi) - M^N_{t^-}(\phi) &= (N\eta_N)^{1/2} \left( \langle \bm{\rho}^N_t, \phi \rangle - \langle \bm{\rho}^N_{t^-}, \phi \rangle \right) \\
	&= (N\eta_N)^{1/2} \left( \langle \rho^N_{Nt/\delta_N^\alpha}, \phi_N \rangle - \langle \rho^N_{(Nt/\delta_N^\alpha)^-}, \phi_N \rangle \right),
	\end{align*}
	using the notation introduced in \eqref{def_phiN}.
	By Definition~\ref{def:slfv}, if $ (t,x_0,r) \in \Pi $,
	\begin{align*}
	\abs{ \langle \rho^N_t, \phi_N \rangle - \langle \rho^N_{t^-}, \phi_N \rangle } &\leq \sup_{k_0 \in [0,1]} \abs{ \int_{\R^d \times [0,1]} \phi_N(x,k) u_N \1_{\abs{x-x_0}<r} (\delta_{k_0}(dk) - \rho^N_{t^-}(x, dk)) dx } \\
	&\leq 2 u_N \int_{\R^d} \delta_N^d \sup_{k \in [0,1]} \abs{ \phi(\delta_N x, k) } \1_{\abs{x-x_0}<r} dx. \numberthis \label{bound_jumps}
	\end{align*}
	By the Cauchy-Schwarz inequality,
	\begin{align*}
	\int_{\R^d} \delta_N^d \sup_{k \in [0,1]} \abs{ \phi(\delta_N x, k) } \1_{\abs{x-x_0}<r} dx &\leq V_r^{1/2} \left( \int_{\R^d} \delta_N^{2d} \sup_{k \in [0,1]} \abs{ \phi(\delta_N x, k) }^2 dx \right)^{1/2} \\
	&\leq \delta_N^{d/2} V_r^{1/2} \| \phi \|_2.
	\end{align*}
	Hence in the fixed radius case,
	\begin{align*}
	\sup_{t \geq 0} \abs{ M^N_t(\phi) - M^N_{t^-}(\phi) } &\leq 2u V_R^{1/2} N^{-1/2} \eta_N^{1/2} \delta_N^{d/2} \|\phi \|_2 \\
	&\leq 2u V_R^{1/2} N^{-1/2} \delta_N \left\| \phi \right\|_2 \cvgas{N} 0.
	\end{align*}
	In the stable case, we use $ \1_{\abs{x-x_0}<r} \leq 1 $ to obtain
	\begin{align*}
	\int_{\R^d} \delta_N^d \sup_{k \in [0,1]} \abs{ \phi(\delta_N x, k) } \1_{\abs{x-x_0}<r} dx \leq \| \phi \|_1
	\end{align*}
	and so
	\begin{align*}
	\sup_{t \geq 0} \abs{ M^N_t(\phi) - M^N_{t^-}(\phi) } \leq 2u N^{-1/2} \| \phi \|_1 \cvgas{N} 0.
	\end{align*}
	This proves \eqref{bound_jumps_MN}.
	For the rest of the proof of Lemma~\ref{lemma:convergence_qvar_MN}, we treat the fixed radius case and the stable case separately.
	
	Recall that
	\begin{align*}
	\left\langle M^N(\phi) \right\rangle_t = u^2 \int_{0}^{t} \langle \eta_N \Gamma^{\nu_\alpha^N}(\bm{\rho}^N_s), \phi \otimes \phi \rangle ds.
	\end{align*}
	In the fixed radius case, this is
	\begin{align} \label{qvar_MN_fixed_radius}
	\left\langle M^N(\phi) \right\rangle_t = u^2 \delta_N^{-2d} \int_{0}^{t} \langle \Gamma^{\delta_N R}(\bm{\rho}^N_s), \phi \otimes \phi \rangle ds,
	\end{align}
	where we define $ \Gamma^{r}(\rho) $ as in \eqref{def_Gamma} with $ \nu = \delta_r $.
	Given the definition of $ [\rho]_r $ in \eqref{def_bracket_rho}, we see that $ \langle \Gamma^{r}(\rho), \phi\otimes \phi \rangle $ contains four terms.
	The first one is
	\begin{multline*}
	\int_{(\R^d)^2} \int_{B(x_1, r) \cap B(x_2, r)} \frac{1}{V_r} \int_{B(y,r)} \int_{[0,1]} \phi(x_1, k) \phi(x_2, k) \rho(z, dk) dz dy dx_1 dx_2 \\ = V_r^2 \int_{\R^d \times [0,1]} \overline{\phi}(y,k,r)^2 \frac{1}{V_r} \int_{B(y,r)} \rho(z,dk) dz dy.
	\end{multline*}
	Introducing the notation
	\begin{align*}
	\overline{\rho}_r(x, dk) := \frac{1}{V_r} \int_{B(x,r)} \rho(z, dk)dz,
	\end{align*}
	this takes the form
	\begin{align*}
	V_r^2 \langle \overline{\rho}_r, \overline{\phi}(\cdot, r)^2 \rangle.
	\end{align*}
	The second term in $ \langle \Gamma^r(\rho), \phi \otimes \phi \rangle $ is
	\begin{multline} \label{Gamma_r_1}
	\int_{(\R^d)^2} \int_{B(x_1, r) \cap B(x_2, r)} \frac{1}{V_r} \int_{B(y,r)} \int_{[0,1]^2} \phi(x_1, k_1) \phi(x_2, k_2) \rho(x_1, dk_1) \rho(z, dk_2) dz dy dx_1 dx_2 \\ = \int_{\R^d} \int_{B(y,r) \times [0,1]} \phi(x_1, k_1) \rho(x_1, dk_1) dx_1 \int_{[0,1]} \int_{B(y,r)} \phi(x_2, k_2) dx_2 \frac{1}{V_r} \int_{B(y,r)} \rho(z, dk_2) dz dy.
	\end{multline}
	Now note that
	\begin{align*}
	\abs{ \int_{B(y,r) \times [0,1]} (\phi(x_1, k_1) - \phi(y, k_1)) \rho(x_1, dk_1) dx_1 } \leq  r V_r \max_{\abs{\beta} = 1} \|\partial_\beta \phi \|_\infty.
	\end{align*}
	Thus, if one replaces $ \phi(x_i, k_i) $ by $ \phi(y, k_i) $ in \eqref{Gamma_r_1}, the difference between the two expressions is at most
	\begin{align*}
	2 r V_r^2 \max_{\abs{\beta} = 1} \|\partial_\beta \phi \|_\infty \| \phi \|_1.
	\end{align*}
	After this substitution, \eqref{Gamma_r_1} becomes
	\begin{align*}
	V_r^2 \int_{\R^d \times [0,1]^2} \phi(y, k_1) \phi(y, k_2) \overline{\rho}_r(y, dk_1) \overline{\rho}_r(y,dk_2) dy =: \langle \overline{\rho}_r \cdot \overline{\rho}_r, \phi \otimes \phi \rangle,
	\end{align*}
	setting
	\begin{align*}
	\overline{\rho}_r \cdot \overline{\rho}_r (y, dk_1 dk_2) := \overline{\rho}_r(y, dk_1) \overline{\rho}_r(y,dk_2).
	\end{align*}
	We also note that the same reasoning applies to the last two terms in $ \langle \Gamma^r(\rho), \phi \otimes \phi \rangle $, modulo the sign in the front.
	It follows that
	\begin{align*}
	\abs{ \langle \Gamma^r(\rho), \phi \otimes \phi \rangle - V_r^2 \left[ \langle \overline{\rho}_r, \phi^2 \rangle - \langle \overline{\rho}_r\cdot \overline{\rho}_r, \phi \otimes \phi \rangle \right] } \leq 8 r V_r^2 \max_{\abs{\beta} = 1} \|\partial_\beta \phi \|_\infty \| \phi \|_1.
	\end{align*}
	Coming back to \eqref{qvar_MN_fixed_radius}, this implies
	\begin{multline*}
	\abs{ \left\langle M^N(\phi) \right\rangle_t - u^2V_R^2 \int_{0}^{t} \left[ \langle \overline{(\bm{\rho}^N_s)}_{\delta_N R}, \phi^2 \rangle - \langle \overline{(\bm{\rho}^N_s)}_{\delta_N R} \cdot \overline{(\bm{\rho}^N_s)}_{\delta_N R}, \phi \otimes \phi \rangle \right] ds } \\ \leq 8 t\, u^2 V_R^2\, R\, \delta_N \max_{\abs{\beta} = 1} \|\partial_\beta \phi \|_\infty\, \| \phi \|_1,
	\end{multline*}
	which vanishes as $ N \to \infty $.
	By Proposition~\ref{prop:convergence_deterministic_limit}, $ (\bm{\rho}^N_t, t \in [0,T]) $ converges in probability to $ \lambda $ as $ N \to \infty $.
	It follows that $ (\overline{(\bm{\rho}^N_t)}_{\delta_N R}, t \geq 0) $ converges to the same limit and that
	\begin{align*}
	\langle \overline{(\bm{\rho}^N_s)}_{\delta_N R} \cdot \overline{(\bm{\rho}^N_s)}_{\delta_N R}, \phi \otimes \phi \rangle \cvgas{N} \int_{\R^d} \left( \int_{[0,1]} \phi(x, k) dk \right)^2 dx
	\end{align*}
	in probability, uniformly for $ s \in [0,t] $.
	As a result, recalling \eqref{def_Q}, for any $ t \geq 0 $,
	\begin{align*}
	\left\langle M^N(\phi) \right\rangle_t \cvgas{N} \langle \mathcal{Q}, \phi \otimes \phi \rangle\, t
	\end{align*}
	in probability.
	Lemma~\ref{lemma:convergence_qvar_MN} is then proved in the fixed radius case.
	
	In the stable case, by the definition of $ \nu_\alpha^N $,
	\begin{align*}
	\left\langle M^N(\phi) \right\rangle_t &= u^2 \int_{0}^{t} \int_{1}^{\infty} \langle \Gamma^{\delta_N r}(\bm{\rho}^N_s), \phi \otimes \phi \rangle \delta_N^{-(\alpha + d)} \frac{dr}{r^{1+\alpha+d}} ds \\
	&= u^2 \int_{0}^{t} \int_{\delta_N}^{\infty} \langle \Gamma^r(\bm{\rho}^N_s), \phi \otimes \phi \rangle \frac{dr}{r^{1+\alpha+d}} ds.
	\end{align*}
	Recall from \eqref{bound_qvar_1} that
	\begin{align*}
	\abs{ \langle \Gamma^r(\rho), \phi \otimes \phi \rangle } &\leq 4 \int_{(\R^d)^2} V_{2,r}(x,y) \sup_{k \in [0,1]} \abs{\phi(x,k)} \sup_{k \in [0,1]} \abs{\phi(y,k)} dx dy \\
	&\leq 4 V_r^2 \left\| \phi \right\|_2^2.
	\end{align*}
	It follows that
	\begin{align*}
	\int_{0}^{t} \int_{0}^{\delta_N} \langle \Gamma^r(\bm{\rho}^N_s), \phi \otimes \phi \rangle \frac{dr}{r^{1+\alpha+d}} ds \cvgas{N} 0,
	\end{align*}
	almost surely.
	We then prove that
	\begin{align} \label{cvg_qvar_stable_2}
	\int_{0}^{t} \int_{0}^{\infty} \langle \Gamma^r(\bm{\rho}^N_s), \phi \otimes \phi \rangle \frac{dr}{r^{1+\alpha+d}} ds \cvgas{N} \int_{0}^{t} \int_{0}^{\infty} \langle \Gamma^r(\lambda), \phi \otimes \phi \rangle \frac{dr}{r^{1+\alpha+d}} ds
	\end{align}
	in probability.
	To do this, we first show that the integrand converges.
	From the definition of $ \Gamma^r(\rho) $, we can write,
	\begin{align*}
	\langle \Gamma^r(\rho), \phi \otimes \phi \rangle &= V_r^2 \langle \rho, \overline{(\overline{\phi})^2}(\cdot, r) \rangle - \langle \rho \otimes \rho, \Psi_r \rangle
	\end{align*}
	where
	\begin{multline*}
	\Psi_r(x_1, x_2, k_1, k_2) = \phi(x_1, k_1) \int_{B(x_1, r) \cap B(x_2, r)} \overline{\phi}(y, k_2, r) dy \\ + \phi(x_2, k_2) \int_{B(x_1, r) \cap B(x_2, r)} \overline{\phi}(y, k_1, r) dy - \phi(x_1, k_2)\phi(x_2, k_2).
	\end{multline*}
	By Proposition~\ref{prop:convergence_deterministic_limit}, $ \bm{\rho}^N_s $ converges to $ \lambda $ in the vague topology in probability and uniformly for $ s \in [0,t] $.
	Hence $ \bm{\rho}^N_s \otimes \bm{\rho}^N_s $ converges in the same sense to $ \lambda \otimes \lambda $ as $ N \to \infty $.
	As a consequence, for any $ r > 0 $,
	\begin{align*}
	\int_{0}^{t}\langle \Gamma^r(\bm{\rho}^N_s), \phi \otimes \phi \rangle ds \cvgas{N} \int_{0}^{t} \langle \Gamma^r(\lambda), \phi \otimes \phi \rangle ds
	\end{align*}
	in probability.
	Then, by the bound on \eqref{bound_qvar_1}, we can use dominated convergence to obtain \eqref{cvg_qvar_stable_2}.
	To conclude, we note that
	\begin{align*}
	\Gamma^{r}(\lambda)(x_1, x_2, dk_1 dk_2) = V_{2,r}(x_1, x_2) \left[ dk_1 \delta_{k_1}(dk_2) - dk_1 dk_2 \right].
	\end{align*}
	Integrating over $ r $, and recalling from \eqref{def_Kalpha} that
	\begin{align*}
		K_\alpha(x,y) = \int_{0}^{\infty} V_{2,r}(x,y) \frac{dr}{r^{d+\alpha+1}},
	\end{align*}
	we see that
	\begin{align*}
		u^2 \int_{0}^{\infty} \langle \Gamma^r(\lambda), \phi \otimes \phi \rangle \frac{dr}{r^{1+\alpha+d}} = \langle \mathcal{Q}_\alpha, \phi \otimes \phi \rangle.
	\end{align*}
	This concludes the proof of Lemma~\ref{lemma:convergence_qvar_MN}.
\end{proof}

	\subsection{Proof of the central limit theorems} \label{subsec:proof_clt}
	
	Let us now prove Theorem~\ref{thm:clt_fixed_radius} and \ref{thm:clt_stable}, applying Theorem~\ref{thm:cvg_S'} to $ (Z^N, N \geq 1) $.
	
	\begin{proof}[Proof of Theorem~\ref{thm:clt_fixed_radius} and \ref{thm:clt_stable}]
		Condition \textit{i}) of Theorem~\ref{thm:cvg_S'} is satisfied by $ Z^N $ by Lemma~\ref{lemma:tightness}.
		To check condition \textit{ii}), let $ \phi_1, \ldots, \phi_p $ be elements of $ \mathcal{S}(\R^d \times [0,1]) $ and $ t_1, \ldots, t_p \in [0,T] $. 
		We shall apply Proposition~\ref{prop:cvg_sto_int}, with
		\begin{align*}
		f^N_i(s,x,k) = \1_{s \leq t_i}\, P^{N,\alpha}_{t_i-s} \phi_i(x,k), && f_i(s,x,k) = \1_{s \leq t_i}\, P^{(\alpha)}_{t_i-s} \phi_i(x,k).
		\end{align*}
		Then, by \eqref{stochastic_integral},
		\begin{align*}
		M^N_T(f^N_i) = \langle Z^N_{t_i}, \phi_i \rangle, && M_T(f_i) = \langle Z_{t_i}, \phi_i \rangle,
		\end{align*}
		where $ (Z_t, t \geq 0) $ is defined by
		\begin{align} \label{int_sto_Z}
		\langle Z_t, \phi \rangle := \int_{[0,t] \times \R^d \times [0,1]} P^{(\alpha)}_{t-s} \phi(x,k) M(ds\, dx\, dk).
		\end{align}
		By Theorem~5.1 in \cite{walsh_introduction_1986}, $ (Z_t, t\geq 0) $ is the unique solution in $ D(\R_+, \mathcal{S}'(\R^d \times [0,1])) $ to the stochastic partial differential equation \eqref{spde_Z_fixed_radius} in the fixed radius case and \eqref{spde_Z_stable} in the stable case.
		The assumptions of Proposition~\ref{prop:cvg_sto_int} are then straightforward to check with the help of Lemma~\ref{lemma:test_functions} and Lemma~\ref{lemma:cvg_MN}.
		This yields the convergence
		\begin{align} \label{cvg_marginals}
		\left( \langle Z^N_{t_1}, \phi_1 \rangle, \ldots, \langle Z^N_{t_p}, \phi_p \rangle \right) \longrightarrow \left( \langle Z_{t_1}, \phi_1 \rangle, \ldots, \langle Z_{t_p}, \phi_p \rangle \right),
		\end{align}
		in distribution as $ N \to \infty $.
		Hence, by Theorem~\ref{thm:cvg_S'}, there exists a process $ (\tilde{Z}_t, t \in [0,T]) $ with sample paths in $ D(\R_+, \mathcal{S}'(\R^d \times [0,1])) $ such that $ Z^N $ converges to $ \tilde{Z} $ in distribution.
		By \eqref{cvg_marginals}, for any $ \phi_1, \ldots, \phi_p $ and $ t_1, \ldots, t_p $ as above, 
		\begin{align*}
		\left( \langle \tilde{Z}_{t_1}, \phi_1 \rangle, \ldots, \langle \tilde{Z}_{t_p}, \phi_p \rangle \right) =_d \left( \langle Z_{t_1}, \phi_1 \rangle, \ldots, \langle Z_{t_p}, \phi_p \rangle \right),
		\end{align*}
		where $ =_d $ stands for equality in distribution.
		It follows that $ \tilde{Z} =_d Z $, and we conclude that $ Z^N $ converges in distribution to $ Z $.
	\end{proof}

	 \subsection{Proof of the central limit theorem in the non-stationary regime with general mutation mechanism} \label{subsec:proof_clt_non-stationary}
	 
	 Let us finish this section by outlining the main adaptations needed to prove Theorem~\ref{thm:clt_non_stationary}.
	 Recall that, between reproduction events,
	 \begin{align*}
	 	\rho^N_t(x,\cdot) = \mathcal{T}^*_{\mu_N (t-s)} \rho^N_s (x, \cdot). 
	 \end{align*}
	 As a result, adapting the proof of Proposition~\ref{prop:slfv_semimartingale}, we obtain
	 \begin{align*}
	 	\langle \rho^N_t, \phi \rangle = \langle \rho^N_0, \phi \rangle + \int_{0}^{t} \left\lbrace \mu_N \langle \rho^N_s, \mathcal{G} \phi \rangle + u_N \int_{0}^{\infty} V_r \langle \rho^N_s, \overline{\overline{\phi}}(\cdot, r) - \phi \rangle \nu_\alpha(dr) \right\rbrace ds + \mathcal{M}^N_t(\phi),
	 \end{align*}
	 and $ (\mathcal{M}^N_t(\phi), t \geq 0) $ is a square-integrable martingale with predictable variation process
	 \begin{align*}
	 	\langle \mathcal{M}^N(\phi) \rangle_t = u_N^2 \int_{0}^{t} \langle \Gamma^{\nu_\alpha} (\rho^N_s), \phi \otimes \phi \rangle ds.
	 \end{align*}
	 Since $ \mathcal{G} $ does not act on the space variable, we have
	 \begin{align*}
	 	\mathcal{G} \phi(\delta_N \cdot, \cdot) = (\mathcal{G} \phi)(\delta_N \cdot, \cdot).
	 \end{align*}
	 As a result, for the rescaled process $ \bm{\rho}^N_t $, using the notation \eqref{def_phiN} from Section~\ref{subsec:semimartingale},
	 \begin{align} \label{rescaled_rho_N_G}
		 \langle \bm{\rho}^N_t, \phi \rangle = \langle \bm{\rho}^N_0, \phi \rangle + \int_0^t \langle \bm{\rho}^N_s, \mu \mathcal{G} \phi + u \mathcal{L}^{N,\alpha} \phi \rangle ds + \mathcal{M}^N_{Nt / \delta_N^\alpha}(\phi_N).
	 \end{align}
	 We now let $ (P^{N,\alpha}_t, t \geq 0) $ be the semigroup acting on continuous and bounded functions on $ \R^d \times [0,1] $ defined by
	 \begin{align*}
		 P^{N,\alpha}_t \phi(x,k) := \E[(x,k)]{\phi(X^{N,\alpha}_t, \mathcal{K}_t)},
	 \end{align*}
	 where $ X^{N,\alpha} $ and $\mathcal{K}$ are two independent Markov processes taking values respectively in $\R^d$ and $[0,1]$ and with generators $ u\mathcal{L}^{N,\alpha} $ and $ \mu \mathcal{G} $, and $\E[(x,k)]{\cdot} $ denotes the expectation with repect to their joint distribution, with the initial condition $ X^{N,\alpha}_0 = x $ and $\mathcal{K}_0 = k $.
	 We then define
	 \begin{align} \label{def_centring_term}
		 \bm{p}^N_t(x,dk) := (P^{N,\alpha}_t)^* \bm{\rho}^N_0 (x, dk).
	 \end{align}
	 With this definition, for any $ \phi \in \mathcal{S}(\R^d \times [0,1]) $,
	 \begin{align*}
		 \langle \bm{p}^N_t, \phi \rangle = \langle \bm{\rho}^N_0, \phi \rangle + \int_0^t \langle \bm{p}^N_s, \mu \mathcal{G} \phi + u \mathcal{L}^{N,\alpha} \phi \rangle ds.
	 \end{align*}
	 Subtracting this to \eqref{rescaled_rho_N_G} and multiplying by $ (N \eta_N)^{1/2} $ yields
	 \begin{align*}
		 \langle Z^N_t, \phi \rangle = \int_0^t \langle Z^N_s, \mu \mathcal{G} \phi + u \mathcal{L}^{N,\alpha} \phi \rangle ds + M^N_t(\phi),
	 \end{align*}
	 where we have set $ M^N_t(\phi) = (N \eta_N)^{1/2} \mathcal{M}^N_{Nt/\delta_N^\alpha}(\phi_N)$.
	 The predictable variation of $ M^N(\phi) $ is unchanged and is given by \eqref{qvar_rescaled}.

	 Lemma~\ref{lemma:bound_qvar_MN} does not require any modification, and we can adapt the proof of Lemma~\ref{lemma:cvg_PN} to obtain the following, where $(P^{(\alpha)}_t, t \geq 0) $ now denotes the semigroup acting on measurable and bounded functions on $\R^d \times [0,1]$ generated by $\mu \mathcal{G} + u \mathcal{D}_\alpha $.

	 \begin{lemma} \label{lemma:cvg_PN_general}
		 For any $ \phi \in \mathcal{S}(\R^d \times [0,1]) $, $t \geq 0$, $N \geq 1$ and $q \geq 1$,
		 \begin{align*}
			 \| P^{N,\alpha}_t \phi \|_q \leq \| \phi \|_q.
		 \end{align*}
		 Moreover, for any $\beta \in \N^d$,
		 \begin{align*}
			 \| \partial_\beta P^{N,\alpha}_t \phi \|_q \leq \| \partial_\beta \phi \|_q,
		 \end{align*}
		 and there exists a constant $\Cst{cvgpsi} > 0 $ such that, for all $N \geq 1$ and $t \geq 0$,
		 \begin{align*}
			 \| P^{N,\alpha}_t \phi - P^{(\alpha)}_t \phi \|_q \leq \Cst{cvgpsi}\, t \, (\delta_N)^\gamma \max_{0 \leq |\beta| \leq 4} \| \partial_\beta \phi \|_q,
		 \end{align*}
		 where $\gamma = 2$ in the fixed radius case and $\gamma=2-\alpha$ in the stable case.
	 \end{lemma}
 
 	With this result, and noting that the definition of $ \bm{p}^{(\alpha)} $ in \eqref{def_p_alpha} is equivalent to
 	\begin{align*}
 		\bm{p}^{(\alpha)}_t(x,dk) = (P^{(\alpha)}_t)^* \bm{\rho}_0 (x,dk),
 	\end{align*}
 	we obtain the following.
 	
 	\begin{proposition} \label{prop:cvg_centring_term}
 		Under the assumptions of Theorem~\ref{thm:clt_non_stationary}, for any $ T > 0 $,
 		\begin{align*}
 			\lim_{N \to \infty} \E{ \sup_{t \in [0,T]} d(\bm{p}^N_t, \bm{p}^{(\alpha)}_t)} = 0.
 		\end{align*}
 		Furthermore, if
 		\begin{align*}
 			\sup_{x \in \R^d} d_{TV}(\bm{\rho}^N_0(x,\cdot), \bm{\rho}_0(x,\cdot)) \leq C (\delta_N)^{\gamma},
 		\end{align*}
 		almost surely for some constant $ C > 0 $ where $ \gamma = 2 $ in the fixed radius case and $ \gamma = 2-\alpha $ in the stable case and $ d_{VT} $ denotes the total variation distance, then for any $ T > 0 $ there exists another constant, still denoted $ C > 0 $, such that, almost surely
 		\begin{align*}
 			\sup_{t \in [0,T]} d(\bm{p}^N_t, \bm{p}_t^{(\alpha)}) \leq C (\delta_N)^{\gamma}.
 		\end{align*}
 	\end{proposition}

	\begin{proof}
		By the definition of the distance $ d $ in \eqref{def:d},
		\begin{align*}
			d(\bm{p}^N_t, \bm{p}^{(\alpha)}_t) = \sum_{n=1}^{\infty} \frac{1}{2^n} \abs{ \langle \bm{p}^N_t, \phi_n \rangle - \langle \bm{p}^{(\alpha)}_t, \phi_n \rangle }.
		\end{align*}
		Then, by the definition of $ \bm{p}^N_t $ in \eqref{def_centring_term} and that of $ \bm{p}^{(\alpha)}_t $ in \eqref{def_p_alpha},
		\begin{align*}
			\abs{\langle \bm{p}^N_t, \phi_n \rangle - \langle \bm{p}^{(\alpha)}_t, \phi_n \rangle} \leq \abs{\langle \bm{\rho}^N_0, P^{N,\alpha}_t \phi_n - P^{(\alpha)}_t \phi_n \rangle} + \abs{\langle \bm{\rho}^N_0, P^{(\alpha)}_t \phi_n \rangle - \langle \bm{\rho}_0, P^{(\alpha)}_t \phi_n \rangle}.
		\end{align*}
		For the first term on the right, we use the definition of the norm $ \| \cdot \|_q $, Lemma~\ref{lemma:cvg_PN_general} and \eqref{phi_n_bound} to obtain
		\begin{align*}
			\abs{\langle \bm{\rho}^N_0, P^{N,\alpha}_t \phi_n - P^{(\alpha)}_t \phi_n \rangle} &\leq \| P^{N,\alpha}_t \phi_n - P^{(\alpha)}_t \phi_n \|_1 \\
			&\leq \Cst{cvgpsi}\, t\, (\delta_N)^\gamma \Cst{phi_n}.
		\end{align*}
		For the second term, we note that
		\begin{align*}
			\abs{\langle \bm{\rho}^N_0, P^{(\alpha)}_t \phi_n \rangle - \langle \bm{\rho}_0, P^{(\alpha)}_t \phi_n \rangle} &\leq 2 \| P^{(\alpha)}_t \phi_n \|_1\\
			&\leq 2 \Cst{phi_n},
		\end{align*}
		using Lemma~\ref{lemma:cvg_PN_general} and \eqref{phi_n_bound}.
		As a result, for any $ \varepsilon > 0 $,
		\begin{multline*}
			\E{ \sup_{t \in [0,T]} d(\bm{p}^N_t, \bm{p}^{(\alpha)}_t) } \leq \Cst{cvgpsi} \Cst{phi_n} \, t\, (\delta_N)^{\gamma} + \varepsilon \\ + 2 \Cst{phi_n} \sum_{n=1}^{\infty} \frac{1}{2^n} \P{\sup_{t \in [0,T]} \abs{\langle \bm{\rho}^N_0, P^{(\alpha)}_t \phi_n \rangle - \langle \bm{\rho}_0, P^{(\alpha)}_t \phi_n \rangle} > \varepsilon}.
		\end{multline*}
		In addition, since $ \bm{\rho}^N_0 $ converges weakly to $ \bm{\rho}_0 $ in probability and $ P^{(\alpha)}_t\phi_n $ is continuous and bounded, for any $ \varepsilon > 0 $, as $ N \to \infty $,
		\begin{align*}
		\P{\sup_{t \in [0,T]} \abs{\langle \bm{\rho}^N_0, P^{(\alpha)}_t \phi_n \rangle - \langle \bm{\rho}_0, P^{(\alpha)}_t \phi_n \rangle} > \varepsilon} \to 0,
		\end{align*}
		for each $ n \geq 1 $.
		Hence, by dominated convergence,
		\begin{align*}
			\limsup_{N \to \infty} \E{ \sup_{t \in [0,T]} d(\bm{p}^N_t, \bm{p}^{(\alpha)}_t) } \leq \varepsilon,
		\end{align*}
		and the conclusion follows by letting $ \varepsilon \to 0 $.
		The second part of the statement is trivial since
		\begin{align*}
			\abs{\langle \bm{\rho}^N_0, P^{(\alpha)}_t \phi_n \rangle - \langle \bm{\rho}_0, P^{(\alpha)}_t \phi_n \rangle} \leq \sup_{x \in \R^d} d_{TV}(\bm{\rho}^N_0(x,\cdot), \bm{\rho}_0(x,\cdot)) \, \| P^{(\alpha)}_t \phi_n \|_1.
		\end{align*}
		This concludes the proof of Proposition~\ref{prop:cvg_centring_term}.
	\end{proof}

	We can then use Theorem~\ref{thm:walsh} as in the proof of Proposition~\ref{prop:convergence_deterministic_limit} to obtain
	\begin{align*}
		\E{ \sup_{t \in [0,T]} d(\bm{\rho}^N_t, \bm{p}^N_t) } \leq \frac{C}{(N \eta_N)^{1/2}},
	\end{align*}
	for some $ C > 0 $, and \eqref{cvg_det_limit_general} follows from the above combined with Proposition~\ref{prop:cvg_centring_term}.
	
	The convergence of the sequence of martingale measures $ (M^N, N \geq 1) $ is identical to the proof of Lemma~\ref{lemma:cvg_MN}, replacing the convergence of $ \bm{\rho}^N_t $ to $ \bm{p}^{(\alpha)}_t $ instead of $ \lambda $ and noting that \eqref{qvar_non_stationary} coincides with
	\begin{align*}
		u^2 \int_{0}^{t} \int_{0}^{\infty} \langle \Gamma^r(\bm{p}^{(\alpha)}_s), \phi \otimes \phi \rangle \frac{dr}{r^{1+d+\alpha}} ds.
	\end{align*}
	The convergence of the sequence $ (Z^N, N \geq 1) $ then follows exactly as in the stationary setting.
	This concludes the proof of Theorem~\ref{thm:clt_non_stationary}.
    
    \section{Derivation of the Wright-Malécot formula} \label{sec:WM_proof}
    
    The aim of this section is to prove Theorem~\ref{thm:WM_formula}.
    Recall that 
    \begin{align*}
    P^N_t(\phi, \psi) = \E{ \left\langle \bm{\rho}^N_t \otimes \bm{\rho}^N_t, (\phi \otimes \psi)\, \1_{\Delta} \right\rangle },
    \end{align*}
    where $ \phi : \R^d \to \R_+ $ and $ \psi : \R^d \to \R_+ $ are two smooth  and compactly supported probability density functions and
    \begin{align*}
    \1_{\Delta}(k_1, k_2) = \1_{k_1 = k_2}.
    \end{align*}
    Since the Lebesgue measure of $ \lbrace (k_1, k_2) \in [0,1]^2 : k_1 = k_2 \rbrace $ is zero, $ \langle \lambda \otimes \lambda, (\phi \otimes \psi)\, \1_\Delta \rangle = 0 $ and this is equivalent to
    \begin{align} \label{PN_ZN}
    P^N_t(\phi, \psi) = (N \eta_N)^{-1} \E{ \left\langle Z^N_t \otimes Z^N_t, (\phi \otimes \psi)\, \1_{\Delta} \right\rangle }.
    \end{align}
    But, by Theorem~\ref{thm:clt_fixed_radius} and Theorem~\ref{thm:clt_stable}, we have the following.
    
    \begin{lemma} \label{lemma:cvg_PN}
    	For any fixed $ t \geq 0 $, and $ \phi $, $ \psi $ satisfying the assumptions of Theorem~\ref{thm:WM_formula},
    	\begin{align} \label{wm_fixed_radius_N_to_inf}
    	\lim_{N \to \infty} (N \eta_N) P^N_t(\phi, \psi) = \E{ \left\langle Z_t \otimes Z_t, (\phi \otimes \psi)\, \1_{\Delta} \right\rangle },
    	\end{align}
    	where $ (Z_t, t \geq 0) $ solves \eqref{spde_Z_fixed_radius} in the fixed radius case and \eqref{spde_Z_stable} in the stable case.
    \end{lemma}
    
    We prove Lemma~\ref{lemma:cvg_PN} below, but first we make the following observation.
    	By the definition of $ Z_t $ and \eqref{int_sto_Z}, for any $ \phi \in \mathcal{S}(\R^d \times [0,1]) $, $ \langle Z_t, \phi \rangle $ is a Gaussian random variable with mean zero and variance
    	\begin{align*}
    	\int_{0}^{t} \langle \mathcal{Q}_\alpha, \left( P^{(\alpha)}_s \phi \right) \otimes \left( P^{(\alpha)}_s \phi \right) \rangle \, ds.
    	\end{align*}
    	Using Lemma~\ref{lemma:bound_qvar_MN} and Lemma~\ref{lemma:test_functions}, we see that, as $ t \to \infty $, this converges to
    	\begin{align} \label{def_Q_infty}
    	\langle \mathcal{Q}_{\alpha}^{(\infty)}, \phi \otimes \phi \rangle := \int_{0}^{+\infty} \langle \mathcal{Q}_\alpha, \left( P^{(\alpha)}_s \phi \right) \otimes \left( P^{(\alpha)}_s \phi \right) \rangle\, ds.
    	\end{align}
    	Hence $ \langle Z_t, \phi \rangle \to \langle Z, \phi \rangle $ in distribution as $ t \to \infty $, where $ Z $ is a Gaussian random field on $ \R^d \times [0,1] $ with covariation measure $ \mathcal{Q}^{(\infty)}_\alpha $.
    	It is straightforward to extend this to vectors of the form $ (\langle Z_{t}, \phi_1 \rangle, \ldots, \langle Z_t, \phi_k \rangle) $ to obtain the following result.
    
    \begin{cor} \label{cor:stationary_Z}
    	Let $ (Z_t, t \geq 0) $ solve either \eqref{spde_Z_fixed_radius} or $ \eqref{spde_Z_stable} $.
    	Then, as $ t \to \infty $, $ Z_t $ converges in distribution in $ \mathcal{S}'(\R^d \times [0,1]) $ to a Gaussian random field, denoted by $ Z $, with covariation measure $ \mathcal{Q}^{(\infty)}_\alpha $, defined by \eqref{def_Q_infty}.
    \end{cor}
    	
    The proof of Theorem~\ref{thm:WM_formula} then goes along the same lines.
    
    \begin{proof}[Proof of Theorem~\ref{thm:WM_formula}]
    	From \eqref{int_sto_Z}, we see that
    	\begin{align*}
    		\E{ \left\langle Z_t \otimes Z_t, (\phi \otimes \psi)\, \1_{\Delta} \right\rangle } = \int_{0}^{t} \langle \mathcal{Q}_\alpha,  \left( P^{(\alpha)}_s \phi \right) \otimes \left( P^{(\alpha)}_s \psi \right)\, \1_\Delta \rangle\, ds,
    	\end{align*}
    	where we have used the fact that, if $ \Psi(x,k) = \phi(x) h(k) $, then $ P^{(\alpha)}_s \Psi(x,k) = (P^{(\alpha)}_s \phi (x)) h(k) $.
    	Since $ |\1_\Delta| \leq 1 $, we can again use Lemma~\ref{lemma:bound_qvar_MN} and Lemma~\ref{lemma:test_functions} to show that the above expression converges to
    	\begin{align*}
    		\langle \mathcal{Q}^{(\infty)}_\alpha, (\phi \otimes \psi)\, \1_\Delta \rangle = \int_{0}^{+\infty} \langle \mathcal{Q}_\alpha, \left( P^{(\alpha)}_s \phi \right) \otimes \left( P^{(\alpha)}_s \psi \right)\, \1_\Delta \rangle\, ds,
    	\end{align*}
    	as $ t \to \infty $.
    	To conclude, we need to show that the right-hand-side coincides with the expressions given in the statement of Theorem~\ref{thm:WM_formula}.
    	
    	Let us first do so in the fixed radius case.
    	Recalling the notation \eqref{def_G_2}, by the definition of $ P^{(\alpha)}_t $, for any $ \phi \in \mathcal{S}(\R^d \times [0,1]) $,
    	\begin{align*}
    		P^{(2)}_t \phi(x,k) = e^{-\mu t} \int_{\R^d} G^{(2)}_{ut}(x-y) \phi(y,k) dy.
    	\end{align*}
    	Also recall that, by \eqref{def_Q},
    	\begin{align*}
    		\mathcal{Q}(dx_1 dk_1 dx_2 dk_2) = d x_1 \delta_{x_1}(dx_2) \left( dk_1 \delta_{k_1}(d k_2) - dk_1 dk_2 \right).
    	\end{align*}
    	As a result, using the convolution rule for Gaussian kernels,
    	\begin{align*}
    		\langle \mathcal{Q}^{(\infty)}_2, (\phi \otimes \psi)\, \1_\Delta \rangle = u^2 V_R^2 \int_{0}^{\infty} \int_{(\R^d)^2} e^{-2\mu s} G_{2us}(x-y) \phi(x) \psi(y) dx dy ds.
    	\end{align*}
    	Then, using \citep[p.~146, Eq.~29]{erdelyi_tables_1954}, for $ \alpha,  p > 0 $,
    	\begin{align*}
    	\int_{0}^{\infty} e^{-pt} t^{\nu - 1} e^{-\alpha / 4t} dt = 2 \left( \frac{\alpha}{4 p} \right)^{\nu/2} K_\nu\left( \sqrt{\alpha p} \right),
    	\end{align*}
    	and this yields \eqref{Wright-Malecot_fixed_radius}.
    	
    	In the stable case, similarly,
    	\begin{align*}
    		P^{(\alpha)}_t \phi (x,k) = e^{-\mu t} \int_{\R^d} G^{(\alpha)}_{ut}(x-y) \phi(y, k) dy,
    	\end{align*}
    	where $ G^{(\alpha)}_t $ was defined in \eqref{def_G_alpha}.
    	Recalling the definition of $ \mathcal{Q}_\alpha $ in \eqref{def_Qalpha}, we obtain
    	\begin{multline*}
    		\langle \mathcal{Q}^{(\infty)}_\alpha, (\phi \otimes \psi) \, \1_\Delta \rangle \\ = u^2 \int_{(\R^d)^4} \int_{0}^{\infty}  e^{-2\mu s}   G^{(\alpha)}_{u s}(x_1-y_1) G^{(\alpha)}_{u s}(x_2-y_2) ds K_\alpha(x_1, x_2) \phi(y_1) \psi(y_2) dx_1 dx_2 dy_1 dy_2.
    	\end{multline*}
    	We then use the fact that, from the $ \alpha $-stability property of $ \mathcal{D}_\alpha $,
    	\begin{align*}
    	G^{(\alpha)}_t(x) = \lambda^{-d/\alpha} G^{(\alpha)}_{t/\lambda}(\lambda^{-1/\alpha} x), \quad \forall \lambda > 0,
    	\end{align*}
    	and simple changes of variables to obtain \eqref{Wright-Malecot_stable}.
    \end{proof}
    
    We now prove Lemma~\ref{lemma:cvg_PN}.
    
    \begin{proof}[Proof of Lemma~\ref{lemma:cvg_PN}]
    	By \eqref{PN_ZN}, \eqref{stochastic_integral} and Proposition~\ref{prop:semimartingale_ZN},
    	\begin{align*}
    	N\eta_N \, P^N_t(\phi, \psi) = \E{ u^2 \int_{0}^{t} \langle \eta_N \Gamma^{\nu^N_\alpha}(\bm{\rho}^N_s), \left( P^{N,\alpha}_{t-s} \phi \right) \otimes \left( P^{N,\alpha}_{t-s} \psi \right)\, \1_\Delta \rangle ds }.
    	\end{align*}
    	Combining Lemma~\ref{lemma:bound_qvar_MN} and Lemma~\ref{lemma:test_functions}, we see that there exists a constant $ C > 0 $ (depending on $ \phi $ and $ \psi $) such that
    	\begin{multline*}
    	\abs{ \langle \eta_N \Gamma^{\nu^N_\alpha}(\bm{\rho}^N_s), \left( P^{N,\alpha}_{t-s} \phi \right) \otimes \left( P^{N,\alpha}_{t-s} \psi \right)\, \1_\Delta \rangle - \langle \eta_N \Gamma^{\nu^N_\alpha}(\bm{\rho}^N_s), \left( P^{(\alpha)}_{t-s} \phi \right) \otimes \left( P^{(\alpha)}_{t-s} \psi \right)\, \1_\Delta \rangle } \\ \leq C (\delta_N)^\gamma (t-s)e^{-2\mu(t-s)}.
    	\end{multline*}
    	Furthermore, we can easily adapt the proof of Lemma~\ref{lemma:cvg_MN} to show that, for any $ t \geq 0 $,
    	\begin{align*}
    	u^2 \int_{0}^{t} \langle \eta_N \Gamma^{\nu^N_\alpha}(\bm{\rho}^N_s), \left( P^{(\alpha)}_{t-s} \phi \right) \otimes \left( P^{(\alpha)}_{t-s} \psi \right)\, \1_\Delta \rangle ds \cvgas{N} \int_{0}^{t} \langle \mathcal{Q}_\alpha, \left( P^{(\alpha)}_s \phi \right) \otimes \left( P^{(\alpha)}_s \psi \right)\, \1_\Delta \rangle ds,
    	\end{align*}
    	in probability.
    	Then using Lemma~\ref{lemma:bound_qvar_MN} and Lemma~\ref{lemma:test_functions} again, we can apply the dominated convergence theorem to show that this convergence holds in expectation.
    	This proves
    	\begin{align*}
    	\lim_{N \to \infty } N\eta_N \, P^N_t(\phi, \psi) = \int_{0}^{t} \langle \mathcal{Q}^\alpha, \left( P^{(\alpha)}_s \phi \right) \otimes \left( P^{(\alpha)}_s \psi \right)\, \1_\Delta \rangle ds.
    	\end{align*}
    	The result then follows from \eqref{int_sto_Z} and the definition of the martingale measure $ M $.
    \end{proof}

    %\printbibliography
    \bibliography{ibd_forwards.bib}

%	 \appendix
	 \begin{appendices}
	 
	 \section[Wiener processes on Rd x 0,1]{Wiener processes on $ \R^d \times [0,1] $} \label{sec:wiener}

	 Recall the definition of $ \mathcal{Q} $ in \eqref{def_Q}:
	 \begin{align*}
	 \mathcal{Q}(dx_1 dk_1 dx_2 dk_2) = u^2 V_R^2 d x_1 \delta_{x_1}(dx_2) \left( dk_1 \delta_{k_1}(d k_2) - dk_1 dk_2 \right).
	 \end{align*}

	 \begin{proposition} \label{prop:wiener_fixed_r}
	 	There exists a Wiener process $ (W(t), t \geq 0) $ taking values in $ D(\R_+$, $\mathcal{S}'(\R^d \times [0,1])) $ such that, for all $ \phi, \psi \in \mathcal{S}(\R^d \times [0,1]) $,
	 	\begin{align*}
	 	\E{ \langle W(t), \phi \rangle \langle W(s), \psi \rangle } = t \wedge s \, \langle \mathcal{Q}, \phi \otimes \psi \rangle.
	 	\end{align*}
	 \end{proposition}
	 
	 \begin{proof}
	 	By the definition of $ \mathcal{Q} $, for any $ \phi \in L^2(\R^d \times [0,1]) $,
	 	\begin{align*}
	 	\langle \mathcal{Q}, \phi \otimes \phi \rangle = u^2 V_R^2 \int_{\R^d \times [0,1]} \left( \phi(x,k) - \int_{[0,1]} \phi(x,k')dk' \right)^2 dx dk \geq 0.
	 	\end{align*}
	 	The corresponding (self-adjoint) covariance operator $ \hat{\mathcal{Q}} $ is then given by
	 	\begin{align*}
	 	\hat{\mathcal{Q}} \phi(x,k) =  u^2 V_R^2 \left( \phi(x,k) - \int_{[0,1]} \phi(x,k')dk' \right),
	 	\end{align*}
	 	and the reproducing kernel of $ (W(t), t \geq 0) $ is given by
	 	\begin{align*}
	 		\left\lbrace \phi \in L^2(\R^d \times [0,1]) : \int_{[0,1]} \phi(x,k) dk = 0 \text{ for all }x \in \R^d \right\rbrace.
	 	\end{align*}
	 	Then, by Proposition~I.4.7 in \cite{prato_stochastic_2014}, if $ (e_i, i \geq 1) $ is a complete orthonormal basis of $ L^{2}(\R^d \times [0,1]) $ and if $ (\beta_{i}, i \geq 1) $ is a sequence of independent standard Brownian motions, then
	 	\begin{align*}
	 		W(t) = \sum_{i \geq 1} \beta_{i}(t)\, \hat{\mathcal{Q}}^{1/2} e_i
	 	\end{align*}
	 	defines a $ \mathcal{S}'(\R^d \times [0,1]) $-valued Wiener process with covariance $ \mathcal{Q} $.
	 	This concludes the proof of Proposition~\ref{prop:wiener_fixed_r}.
	 \end{proof}
 
 	 Note that it is easy to construct a complete orthonormal basis $ (e_i, i \geq 1) $ on which $ \hat{\mathcal{Q}} $ is diagonal, simply by taking it of the form $ e_{i,j} = f_i \otimes g_j $ where $ (f_i, i \geq 1) $ (resp. $ (g_j, j \geq 1) $) is a complete orthonormal basis of $ L^{2}(\R^d) $ (resp. of $ L^{2}([0,1]) $) and with $ g_1 = 1 $.
	 
	 Proposition~\ref{prop:wiener_fixed_r} is extended in a straightforward manner to Wiener process with spatial correlations as follows.
	 Recall from \eqref{def_Qalpha} that, for $ \alpha \in (0, 2\wedge d) $,
	 \begin{align*}
	 \mathcal{Q}_\alpha (dx_1 dk_1 dx_2 dk_2) = u^2 \frac{C_{d,\alpha}}{\|x_1 - x_2 \|^{\alpha}} dx_1 dx_2 ( dk_1 \delta_{k_1}(dk_2) - dk_1 dk_2).
	 \end{align*}
	 
	 \begin{proposition} \label{prop:wiener_stable}
	 	For any $ \alpha \in (0,2\wedge d) $, there exists a Wiener process $ (W(t), t \geq 0) $ taking values in $ D([0,T], \mathcal{S}'(\R^d \times [0,1])) $ such that, for all $ \phi, \psi \in \mathcal{S}(\R^d \times [0,1]) $,
	 	\begin{align*}
	 	\E{\langle W(t), \phi \rangle \langle W(s), \psi \rangle} = t \wedge s \, \langle \mathcal{Q}_\alpha, \phi \otimes \psi \rangle.
	 	\end{align*}
	 \end{proposition}
	 
	 The proof is identical to that of Proposition~\ref{prop:wiener_fixed_r} and is omitted.
	 
	 \section{Approximating the Laplacian and the fractional Laplacian} \label{sec:approx_Laplace}
	 
	 The following was proved in \citep[Proposition~A.2]{forien_central_2017} for functions defined on $ \R^d $ and is easily generalised to functions defined on $ \R^d \times [0,1] $.
	 
	 \begin{proposition} \label{prop:averages}
	 	Let $ \phi : \R^d \times [0,1] \to \R $ be twice continuously differentiable with respect to the space variable and suppose that $ \| \partial_\beta \phi \|_q < \infty $ for $ 0 \leq \abs{\beta} \leq 2 $ for some $ 1 \leq q \leq \infty $.
	 	Then, for all $ r > 0 $,
	 	\begin{align*}
	 		\left\| \overline{\phi}(\cdot, r) - \phi \right\|_q \leq \frac{d}{2} r^2 \max_{\abs{\beta} = 2} \left\| \partial_\beta \phi \right\|_q.
	 	\end{align*}
	 	If in addition $ \phi $ admits continuous and $ \| \cdot \|_q $-bounded spatial derivatives of order up to four,
	 	\begin{align} \label{convergence_LN_fixed_radius}
	 		\left\| \overline{\overline{\phi}}(\cdot, r) - \phi - \frac{r^2}{d+2} \Delta \phi \right\|_q \leq \frac{d^3}{3} r^4 \max_{\abs{\beta} = 4} \left\| \partial_\beta \phi \right\|_q.
	 	\end{align}
	 \end{proposition}
 	 
 	 In addition, the following was proved in \citep[Proposition~A.3]{forien_central_2017} for functions defined on $ \R^d $ and is likewise easily generalised to functions defined on $ \R^d \times [0,1] $ (in \cite{forien_central_2017}, the result is stated with $ q \in \lbrace 1, \infty \rbrace $ but the proof also applies to $ 1 \leq q \leq \infty $).
 	 
 	 \begin{proposition} \label{prop:fractional_laplacian}
 	 	For $ \alpha \in (0,d\wedge 2) $, let $ \mathcal{L}^{N,\alpha} $ be the operator defined in \eqref{def_LN}.
 	 	Let $ \phi : \R^d \times [0,1] \to \R $ be twice continuously differentiable and suppose that $ \left\| \partial_\beta \phi \right\|_q < \infty $ for all $ 0 \leq \abs{\beta} \leq 2 $ for some $ 1 \leq q \leq \infty $.
 	 	Then, for all $ \alpha \in (0, 2 \wedge d) $, there exists constants $ \newCst{fracL}, \newCst{approxfracL} > 0 $ which do not depend on $ \phi $ such that, for all $ N \geq 1 $,
 	 	\begin{align*}
 	 		\left\| \mathcal{L}^{N,\alpha} \phi \right\|_q \leq \Cst{fracL} \left( \left\| \phi \right\|_q + \max_{\abs{\beta} = 2} \left\| \partial_\beta \phi \right\|_q \right),
 	 	\end{align*}
 	 	and
 	 	\begin{align} \label{convergence_LN_stable}
 	 		\left\| \mathcal{L}^{N,\alpha}\phi - \mathcal{D}_\alpha \phi \right\|_q \leq \Cst{approxfracL} (\delta_N)^{2-\alpha} \max_{\abs{\beta} = 2} \left\| \partial_\beta \phi \right\|_q.
 	 	\end{align}
 	 \end{proposition}
  
  	 \section{Martingale convergence theorem} \label{sec:Jacod-Sh}
  	 
  	 Here we recall the following result, which can be found in \cite{jacod_limit_2003}.
  	 
  	 \begin{theorem}[Theorem~VIII~3.11 in \cite{jacod_limit_2003}] \label{thm:Jacod_Sh}
  	 	Suppose that $ (X_t, t \geq 0) $, $ X_t = (X^1_t, \ldots, X^d_t) $ is a continuous $ d $-dimensional Gaussian martingale and that for each $ n \geq 1 $, $ (X^n_t, t\geq 0) $ is a \cadlag, locally square-integrable $ d $-dimensional martingale such that:
  	 	\begin{enumerate}[i)]
  	 		\item $ \sup_{t \geq 0} \abs{X^n_t - X^n_{t^-}} $ is bounded uniformly for $ n \geq 1 $ and converges in probability to 0 as $ n \to \infty $,
  	 		\item for each $ t \in \Q $, $ \langle X^{n,i}, X^{n,j} \rangle_t \cvgas{n} \langle X^i, X^j \rangle_t $ in probability.
  	 	\end{enumerate}
   	 Then $ X^n $ converges to $ X $ in distribution in $ D(\R_+, \R^d) $.
  	 \end{theorem}
   
     \section{Mitoma's theorem} \label{sec:walsh}
     
     The following result is due to Mitoma \citep[Theorem~4.1]{mitoma_tightness_1983}, and can also be found in \citep[Theorem~6.13]{walsh_introduction_1986}.
     
     \begin{theorem}[Mitoma's theorem] \label{thm:mitoma}
     	Let $ (X^N_t, t \in [0,T]) $ be a sequence of processes whose sample paths are in $ D([0,T], \mathcal{S}'(\R^d \times [0,1])) $ almost surely.
     	Then the sequence $ \lbrace X^N, N \geq 1 \rbrace $ is tight if and only if, for each $ \phi \in \mathcal{S}(\R^d \times [0,1]) $, the sequence of real-valued processes $ \lbrace (\langle X^N_t, \phi \rangle, t \in [0,T]), N \geq 1 \rbrace $ is tight in $ D([0,T], \R) $.
     \end{theorem}

  	 \section{Convergence of convolution integrals} \label{sec:proof_walsh}
  	 
  	 Let us first give the proof of Theorem~\ref{thm:walsh}, which is adapted from that of Theorem~7.13 in \cite{walsh_introduction_1986}.
  	 
  	 \begin{proof}[Proof of Theorem~\ref{thm:walsh}]
  	 	We first extend $ \psi^N $ to $ \R_+ \times \R_+ \times \R^d \times [0,1] $ by setting, for $ s > t $,
  	 	\begin{align*}
  	 		\psi^N_{s,t} = \psi^N_{t,t}.
  	 	\end{align*}
  	 	Then the extended $ \psi^N $ still satisfies \textit{i)}-\textit{iii)} with the obvious modifications (in particular replacing $ e^{-\mu(t-s)} $ by $ e^{-\mu\max(t-s, 0)} $ in the right-hand-sides in \eqref{condition_psi_1} and \eqref{condition_psi_2}).
  	 	Then fix $ T > 0 $ and define, for $ t \in [0,T] $,
  	 	\begin{align*}
  	 		V^N_t = \int_{[0,T] \times \R^d \times [0,1]} \psi_{s,t}^N(x,k) M^N(ds dx dk).
  	 	\end{align*}
  	 	Then, by \eqref{condition_K_N} and \eqref{condition_psi_2}, for $ t, t' \in [0,T] $,
  	 	\begin{align*}
  	 		\E{ \abs{V^N_{t'} - V^N_t}^2 } &\leq \Cst{KN} \int_{0}^{T} \sum_{q=1}^{k} \left\| \psi_{s,t'}^N - \psi_{s,t}^N \right\|_q^2 ds \\
  	 		&\leq k \Cst{KN} \Cst{psiNtime}^2 \int_{0}^{T} e^{-2\mu (t \wedge t' -s)^+} \abs{t'-t}^2 ds \\
  	 		&\leq k \Cst{KN} \Cst{psiNtime}^2 T \abs{t'-t}^2.
  	 	\end{align*}
  	 	Then, by Kolmogorov's continuity criterion \citep[Corollary~1.2]{walsh_introduction_1986}, for all $ N \geq 1 $ and for any $ \beta \in (0,1/2) $, there exists a random variable $ Y_N > 0 $ such that, for all $ t, t' \in [0,T] $,
  	 	\begin{align} \label{holder_V_N}
  	 		\abs{ V^N_{t'} - V^N_t } \leq Y_N \abs{t'-t}^\beta \qquad \text{ almost surely}
  	 	\end{align}
  	 	and, for all $ N \geq 1 $,
  	 	\begin{align*}
  	 		\E{ Y_N^2 } \leq \newCst{Holder} k \Cst{KN} \Cst{psiNtime}^2 T,
  	 	\end{align*}
  	 	for some constant $ \Cst{Holder} > 0 $.
  	 	
  	 	Let $ (\F^N_t, t \geq 0) $ denote the natural filtration associated to the martingale measure $ M^N $.
  	 	Then
  	 	\begin{align*}
  	 		U^N_t = \E{ V^N_t }{\F^N_t}.
  	 	\end{align*}
  	 	It follows that, for all $ t \in [0,T] $,
  	 	\begin{align*}
  	 		\abs{U^N_t} \leq \E{ \sup_{s \in [0,T]} \abs{V^N_s} }{\F^N_t}.
  	 	\end{align*}
  	 	Noting that the right hand side is a local martingale, we can apply Doob's maximal inequality to write
  	 	\begin{align*}
  	 		\E{ \sup_{t \in [0,T]} \abs{U^N_t}^2 } &\leq 4 \E{ \E{ \sup_{s \in [0,T]} \abs{V^N_s} }{ \F^N_T}^2 } \\
  	 		&\leq 4 \E{ \sup_{t \in [0,T]} \abs{V^N_t}^2 }. \numberthis \label{Doob}
  	 	\end{align*}
  	 	We then use \eqref{holder_V_N} to obtain
  	 	\begin{align*}
  	 		\sup_{t \in [0,T]} \abs{V^N_t} \leq Y_N T^\beta + \abs{V^N_0}.
  	 	\end{align*}
  	 	But, from \eqref{condition_K_N} and \eqref{condition_psi_1},
  	 	\begin{align*}
  	 		\E{ \abs{V^N_0}^2 } &\leq \Cst{KN} \int_{0}^{T} \sum_{q=1}^{k} \left\| \psi^N_{s,0} \right\|_q^2 ds \\
  	 		&\leq \Cst{KN} T k \Cst{psiN}^2.
  	 	\end{align*}
  	 	From which it follows (together with \eqref{holder_V_N}) that
  	 	\begin{align*}
  	 		\E{ \sup_{t \in [0,T]} \abs{V^N_t}^2 }^{1/2} \leq \left( k \Cst{KN} T \right)^{1/2} \left( (\Cst{Holder})^{1/2} \Cst{psiNtime} T^\beta + \Cst{psiN} \right).
  	 	\end{align*}
  	 	Together with \eqref{Doob}, this concludes the proof of \eqref{sup_bound_U}.
  	 	
  	 	It remains to prove the tightness of the sequence $ (U^N, N \geq 1) $ in $ D(\R_+, \R) $.
  	 	To to this, we use Aldous' criterion \cite{aldous_stopping_1978} as stated for processes indexed by $ \R_+ $ in \citep[Theorem~VI~4.5]{jacod_limit_2003}.
  	 	More precisely, we check that the following two conditions are satisfied.
  	 	For $ T > 0 $, let $ \mathcal{T}^N_T $ denote the set of all $ \mathcal{F}^N $-stopping times that are bounded by $ T $.
  	 	\begin{enumerate}[i)]
  	 		\item For all $ T \in \N^* $, $ \varepsilon > 0 $, there exist $ N_0 \in \N^* $, $ K \in \R_+ $ such that
  	 		\begin{align*}
  	 			N \geq N_0 \implies \P{ \sup_{t \in [0,T]} \abs{U^N_t} > K } \leq \varepsilon.
  	 		\end{align*}
  	 		\item For all $ T \in \N^* $, $ \varepsilon > 0 $,
  	 		\begin{align*}
  	 			\lim_{\theta \downarrow 0} \limsup_{N \to \infty} \sup_{S_1, S_2 \in \mathcal{T}^N_T : S_1 \leq S_2 \leq S_1 + \theta} \P{ \abs{U^N_{S_2} - U^N_{S_1}} > \varepsilon } = 0.
  	 		\end{align*}
  	 	\end{enumerate}
   		Condition \textit{i)} clearly follows from \eqref{sup_bound_U} and the Markov inequality.
   		To prove \textit{ii)}, let $ S_1 $ and $ S_2 $ be two stopping times in $ \mathcal{T}^N_{T} $ such that $ S_1 \leq S_2 \leq S_1 + \theta $.
   		By the optional sampling theorem,
   		\begin{align*}
   			U^N_{S_2} - U^N_{S_1} &= \E{ V^N_{S_2} }{ \F^N_{S_2} } - \E{ V^N_{S_1} }{ \F^N_{S_1} } \\
   			&= \E{ V^N_{S_2} - V^N_{S_1} }{ \F^N_{S_2} } + \E{ V^N_{S_1} }{ \F^N_{S_2} } - \E{ V^N_{S_1} }{ \F^N_{S_1} }. \numberthis \label{decomposition_U}
   		\end{align*}
   		By \eqref{holder_V_N},
   		\begin{align*}
   			\abs{ V^N_{S_2} - V^N_{S_1} } \leq Y_N \theta^\beta.
   		\end{align*}
   		On the other hand, by the definition of $ V^N $,
   		\begin{align*}
   			\E{ V^N_{S_1} }{ \F^N_{S_2} } - \E{ V^N_{S_1} }{ \F^N_{S_1} } = \int_{[S_1, S_2] \times \R^d \times [0,1]} \psi^N_{s,S_1}(x,k) M^N(ds dx dk).
   		\end{align*}
   		Note that $ s > S_1 $ in the integral, so the integrand is always $ \psi^N_{S_1, S_1}(x,k) $.
   		Using \eqref{condition_K_N} and \eqref{condition_psi_1}, we then have
   		\begin{align*}
   			\E{ \abs{\E{ V^N_{S_1} }{ \F^N_{S_2} } - \E{ V^N_{S_1} }{ \F^N_{S_1} }}^2 } &\leq \Cst{KN} \theta \E{ \sum_{q=1}^{k} \left\| \psi^N_{S_1, S_1} \right\|_q^2 } \\
   			& \leq \Cst{KN} \Cst{psiN}^2 k \theta.
   		\end{align*}
   		Coming back to \eqref{decomposition_U}, we have shown that there exists a constant $ \newCst{tightness} > 0 $ such that,  for all $ N \geq 1 $,
   		\begin{align*}
   			\E{ \abs{U^N_{S_2} - U^N_{S_1}}^2 }^{1/2} \leq \Cst{tightness} \left( \theta^{\beta} + \theta^{1/2} \right).
   		\end{align*}
   		And \textit{ii)} follows, concluding the proof of Theorem~\ref{thm:walsh}.
  	 \end{proof}
   
   The following result is adapted in the same way from Proposition~7.12 in \cite{walsh_introduction_1986}, and we omit its proof.
   
   \begin{proposition} \label{prop:cvg_sto_int}
   	Let $ (M^N, N \geq 1) $ be a sequence of worthy martingale measures on $ \R^d \times [0,1] $ with dominating measures $ (K_N, N \geq 1) $, and which converges in distribution to $ M $ in $ D([0,T], \mathcal{S}'(\R^d \times [0,1])) $.
   	Suppose that the dominating measures $ K_N $ satisfy \eqref{condition_K_N} for some constant $ \Cst{KN} $ which does not depend on $ N $.
   	Let $ (f^N_i, 1 \leq i \leq p, N \geq 1) $ be a collection of deterministic real-valued functions on $ [0,T] \times \R^d \times [0,1] $ such that
   	\begin{enumerate}[i)]
   		\item \label{fN_schwartz} for each $ N \geq 1 $, $ 1 \leq i \leq p $ and $ s \in [0,T] $, $ f^N_i(s,\cdot) \in \mathcal{S}(\R^d \times [0,1]) $,
   		\item \label{fN_bound_Lq} $ \sup_{N \geq 1} \sup_{s \in [0,T]} \| f^N_i(s,\cdot) \|_q < +\infty $, for $ q \in [1,k] $ and $ 1 \leq i \leq p $,
   		\item \label{fN_cvg} there exist $ f_i : [0,T] \times \R^d \times [0,1] \to \R $, $ 1 \leq i \leq p $, such that, for any $ q \in [1,k] $,
   		\begin{align*}
   			\lim_{N \to \infty} \sup_{s \in [0,T]} \| f^N_i(s,\cdot) - f_i(s,\cdot) \|_q = 0.
   		\end{align*}
   	\end{enumerate}
    Then, letting $ M^N_t(f^N_i) = \int_{[0,t] \times \R^d \times [0,1]} f^N_i(s,x,k) M^N(ds\, dx\, dk) $,
    \begin{align*}
    	\left( M^N_t(f^N_1), \ldots, M^N_t(f^N_p) \right)_{t \in [0,T]} \cvgas{N} \left( M_t(f_1), \ldots, M_t(f_p) \right)_{t \in [0,T]}
    \end{align*}
    in distribution in $ D([0,T], \R^p) $, where $ M_t(f_i) = \int_{[0,t] \times \R^d \times [0,1]} f_i(s, x, k) M(ds\, dx\, dk) $.
    \end{proposition}

    \begin{remark*}
    	In \citep[Proposition~7.12]{walsh_introduction_1986}, the convergence of $ M^N(f^N_i) $ is stated in $ D([0,T]$, $\mathcal{S}'(\R^d)) $, and the $ f^N_i $ are not assumed to take values in $ \mathcal{S}(\R^d) $.
    	Since the statement is proved by integrating against a test function in $ \mathcal{S}(\R^d) $, the content of the statement is essentially identical.
    \end{remark*}

\end{appendices}

\end{document}